\documentclass[12pt,a4paper]{article}

\usepackage[margin = 0.8in]{geometry}
\usepackage{amsfonts,graphicx,subcaption,amsmath,amssymb,amsthm,bbm,hyperref,nicefrac,enumitem,comment}
\usepackage{xcolor}
\hypersetup{
	colorlinks=true
}

\usepackage{cite}
\usepackage{cleveref} % clever reference
\crefname{figure}{Figure}{Figures}
\crefname{lemma}{Lemma}{Lemma}
\crefname{remark}{Remark}{Remark}
\crefname{setting}{Framework}{Framework}
\crefname{theorem}{Theorem}{Theorem}
\crefname{prop}{Proposition}{Proposition}
\crefname{corollary}{Corollary}{Corollary}
\crefname{section}{Section}{Section}
\crefname{subsection}{Subsection}{Subsections}
\crefname{table}{Table}{Tables}

\newcommand{\Cplusplus}{\mbox{C\hspace{-.05em}\raisebox{.4ex}{\tiny\bf +}\hspace{-.10em}\raisebox{.4ex}{\tiny\bf +}}}

\usepackage{listings}
\usepackage{color}
\definecolor{mygreen}{rgb}{0,0.6,0}
\definecolor{mygray}{rgb}{0.5,0.5,0.5}
\definecolor{mymauve}{rgb}{0.58,0,0.82}

\usepackage[super]{nth}

\usepackage{xcolor}
\usepackage[textsize=tiny,textwidth=0.7in,         ]{todonotes}

\usepackage[nolabel]{showlabels}

\newtheorem{lemma}{Lemma}

\newtheorem{setting}[lemma]{Framework}

\usepackage{chngcntr}
\counterwithin*{section}{part}

\newcounter{tempsection}
\newcommand{\sectionbypart}{%
  \setcounter{tempsection}{\value{section}}% Store section counter
  \setcounter{section}{0}% Reset section counter
  \renewcommand{\thesection}{      \arabic{section}}% Section counter display
}

\renewcommand*\thesection{\arabic{section}}
\usepackage{makecell}

\newcommand{\R}{\mathbb{R}}
\newcommand{\N}{\mathbb{N}}
\newcommand{\E}{\mathbb{E}}
\renewcommand{\P}{\mathbb{P}}

\newcommand{\Z}{\mathbb{Z}}

\newcommand{\smallsum}{\textstyle\sum}

\newcommand{\norm}[1]{ \left\| #1 \right\| }

\newcommand{\normmm}[1]{{\left\vert\kern-0.25ex\left\vert\kern-0.25ex\left\vert #1 
    \right\vert\kern-0.25ex\right\vert\kern-0.25ex\right\vert}} %norm with tripple line
\newcommand{\del}{\partial}
\newcommand{\Trace}{\operatorname{Trace}}
\newcommand{\Hess}{\operatorname{Hess}}
\newcommand{\diag}{\operatorname{diag}}

\begin{document}

\title{Numerical simulations for full history recursive \\ multilevel Picard approximations for systems of \\high-dimensional partial differential equations}

\author{
	Sebastian Becker$^1$,
	Ramon Braunwarth$^2$,
	Martin Hutzenthaler$^3$, \\
	Arnulf Jentzen$^{4, 5}$,
	Philippe von Wurstemberger$^6$
	\bigskip
	\\
	\small{$^1$Risklab, Department of Mathematics, ETH Zurich, }  
	\\
	\small{8092 Z\"urich, Switzerland, e-mail: sebastian.becker@math.ethz.ch}
	\smallskip
	\\
	\small{$^2$Department of Mathematics, ETH Zurich,}
	\\
	\small{8092 Z\"urich, Switzerland, e-mail: r.braunwarth@bluewin.ch}
	\smallskip
	\\ 
	\small{$^3$Faculty of Mathematics, University of Duisburg-Essen,}
	\\
	\small{45117 Essen, Germany, e-mail: martin.hutzenthaler@uni-due.de}
	\smallskip
	\\
	\small{$^4$SAM, Department of Mathematics, ETH Zurich,}
	\\
	\small{8092 Z\"urich, Switzerland, e-mail: arnulf.jentzen@sam.math.ethz.ch} 
	\smallskip
	\\
	\small{$^5$Faculty of Mathematics and Computer Science, University of M\"unster,}
	\\
	\small{48149 M\"unster, Germany, e-mail: ajentzen@uni-muenster.de} 
	\smallskip
	\\
	\small{$^6$Risklab, Department of Mathematics, ETH Zurich,}
	\\
	\small{8092 Z\"urich, Switzerland, e-mail: philippe.vonwurstemberger@math.ethz.ch}	
}

\maketitle

\begin{abstract}
One of the most challenging issues in applied mathematics is to develop and analyze algorithms which are able to approximately compute solutions of high-dimensional nonlinear partial differential equations (PDEs).
In particular, it is very hard to develop approximation algorithms which do not suffer under the curse of dimensionality in the sense that the number of computational operations needed by the algorithm to compute an approximation of accuracy $\varepsilon > 0$ grows at most polynomially in both the reciprocal $\nicefrac{1}{\varepsilon}$ of the required accuracy and the dimension $d \in \N$ of the PDE.
Recently, a new approximation method, the so-called \emph{full history recursive multilevel Picard} (MLP) approximation method, has been introduced and, until today, this approximation scheme is the only approximation method in the scientific literature which has been proven to overcome the curse of dimensionality in the numerical approximation of semilinear PDEs with general time horizons.
It is a key contribution of this article to extend the MLP approximation method to systems of semilinear PDEs and to numerically test it on several example PDEs.
More specifically, we apply the proposed MLP approximation method in the case of 
	Allen-Cahn PDEs,
	Sine-Gordon-type PDEs,
	systems of coupled semilinear heat PDEs, and
	semilinear Black-Scholes PDEs
in up to 1000 dimensions.
The presented numerical simulation results suggest in the case of each of these example PDEs that the proposed MLP approximation method produces very accurate results in short runtimes and, in particular, the presented numerical simulation results indicate that the proposed MLP approximation scheme significantly outperforms certain deep learning based approximation methods for high-dimensional semilinear PDEs.
\end{abstract}

\tableofcontents
\sectionbypart

%%%%%%%%%%%%%%%%%%%%%%%%%%%%%%%%%%%%%%%%%%%
%%%%%%%%%%%%%%%%%%%%%%%%%%%%%%%%%%%%%%%%%%%
\section{Introduction} 
%%%%%%%%%%%%%%%%%%%%%%%%%%%%%%%%%%%%%%%%%%%
%%%%%%%%%%%%%%%%%%%%%%%%%%%%%%%%%%%%%%%%%%%

One of the most challenging issues in applied mathematics is to develop and analyze algorithms which are able to approximately compute solutions of high-dimensional nonlinear partial differential equations (PDEs).
In particular, it is very hard to develop approximation algorithms which do not suffer under the curse of dimensionality in the sense that the number of computational operations needed by the algorithm to compute an approximation of accuracy $\varepsilon > 0$ grows at most polynomially in both the reciprocal $\nicefrac{1}{\varepsilon}$ of the required accuracy and the dimension $d \in \N$ of the PDE.
In the last four years, very significant progress has been made in this research area, where particularly the following two types of approximation methods have turned out to be very promising:

\begin{enumerate}[label=(\Roman{*})]
	
	\item \label{intro_approaches:item1}
	Deep learning based approximation methods for PDEs; cf., e.g., 
	\cite{Becker2018,Becker2019solving,beck2017machine,BeckBeckerGrohsJaafariJentzen18,beck2019deep,berg2018unified,chan2018machine,chen2019deep,chiaramonte2013solving,dockhorn2019discussion,Farahmand2017deep,Fujii19Asymptotic,Goudenege2019Machine,Han2018,han2018solving,he2018relu,henry2017deep,hure2019some,Jacquier2019,jianyu2003numerical,khoo2017solving,lagaris1998artificial,lee1990neural,Long2017,lye2020deep,magill2018neural,meade1994numerical,nabian2018deep,Nusken2020solving,pham2019neural,raissi2018deep,raissi2018forward,ramuhalli2005finite,sirignano2018dgm,uchiyama1993solving,E17Ritz,Weinan2017,becker2019pricing}
	
	\item \label{intro_approaches:item2}
	Full history recursive multilevel Picard approximation methods for PDEs; cf., e.g., 
	\cite{EHutzenthalerJentzenKruse16,hutzenthaler2017multi,OvercomingPaper,hutzenthaler19gradientdependent,E2019,Beck2020MLPAC,Giles19Generalized,VW2019defaultRisk,Beck2020MLPelliptic}
	(in the following we abbreviate \emph{full history recursive multilevel Picard} by MLP)
\end{enumerate}
Roughly speaking, deep learning based approximation methods for high-dimensional PDEs 
% have been introduced in \cite{han2018solving,Weinan2017} and 
 are often based on the idea

\begin{enumerate}[label=(\roman{*})]
	\item[(Ia)] \label{DL_idea:item1}
	to approximate the solution of the considered PDE through the solution of a suitable infinite dimensional stochastic optimization problem on an appropriate function space,
	
	\item[(Ib)] \label{DL_idea:item2}
	to  approximate some of the functions appearing in the infinite dimensional stochastic optimization problem by deep neural networks (DNNs) to obtain finite dimensional stochastic optimization problems, and
	
	\item[(Ic)] \label{DL_idea:item3}
	to apply stochastic gradient descent type algorithms to the resulting finite dimensional stochastic optimization problems to approximately learn the optimal parameters of the involved DNNs.
	
\end{enumerate}
MLP approximation methods have first been proposed in \cite{EHutzenthalerJentzenKruse16,OvercomingPaper} and are, roughly speaking, based on the idea

\begin{enumerate}[label=(\roman{*})]
	\item[(IIa)] \label{MLP_idea:item1}
	to reformulate the computational problem under consideration as a stochastic fixed point equation on a suitable function space with the fixed point of the fixed point equation being the solution of the computational problem,
	
	\item[(IIb)] \label{MLP_idea:item2}
	to approximate the solution of the fixed point equation by means of iterations given by the fixed point equation (which are referred to as Picard iterations in the context of temporal integral fixed point equations), and
	
	\item[(IIc)] \label{MLP_idea:item3}
	to approximate the resulting fixed point iterates by suitable multilevel Monte Carlo approximations, which are full history recursive in the sense that for all $n \in \N$ it holds that the multilevel Monte Carlo approximation of the $n$th fixed point iterate is based on evaluations of multilevel Monte Carlo approximations of the $(n-1)$th, $(n-2)$th, \ldots, 2nd, and 1st fixed point iterates.
\end{enumerate}
A key advantage of deep learning based approximation methods for PDEs is that they seem to be applicable to a very wide class of PDEs including 
	semilinear parabolic PDEs (cf, e.g., \cite{beck2019deep,han2018solving,Weinan2017}),
	elliptic PDEs (cf, e.g., \cite{jianyu2003numerical,E17Ritz}),
	free boundary PDEs associated to optimal stopping problems (cf, e.g., \cite{Becker2018,Becker2019solving,becker2019pricing,Goudenege2019Machine,chen2019deep}), and
	fully nonlinear PDEs (cf, e.g., \cite{beck2017machine,pham2019neural}),
while MLP approximation algorithms are limited to the situation where the computational problem can be formulated as a suitable stochastic fixed point equation and thereby (currently) exclude, for example, fully nonlinear PDEs.
On the other hand, a key advantage of MLP approximation methods is that, until today, these approximation methods are the only methods for which it has been proven that they overcome the curse of dimensionality in the numerical approximation of semilinear PDEs with general time horizons.
In contrast to this, for deep learning based approximation methods there are so far a number of encouraging numerical simulations for PDEs, but only partial error analysis results 
(see, e.g., \cite{GrohsHornungJentzenVW18,Elbraechter2018DNN,Jentzen2018Aproof,Berner2018Analysis,Han2018,Gonon2019Uniform,Grohs2019Spacetime,Grohs2019Deep,Hutzenthaler2019Aproof,Kutyniok2019Atheoretical,Reisinger2019Rectified}) 
which corroborate the conjecture that deep learning based approximation methods might overcome the curse of dimensionality.
These partial error analysis results prove that there exist DNNs which are able to approximate solutions of PDEs with the number of parameters in the DNN growing at most polynomially in both the PDE dimension and the reciprocal of the required approximation accuracy, but there are no results asserting that the employed stochastic optimization algorithm will find such an approximating DNN.
Moreover, one should note that in the case of nonlinear PDEs the proofs for the partial error analysis results mentioned above (cf.\ \cite{Hutzenthaler2019Aproof}) are, in turn, strongly based on the fact that MLP approximation schemes overcome the curse of dimensionality (cf.\ \cite{OvercomingPaper}).

The above mentioned articles \cite{EHutzenthalerJentzenKruse16,hutzenthaler2017multi,OvercomingPaper,hutzenthaler19gradientdependent,Beck2020MLPAC,Giles19Generalized,VW2019defaultRisk,Beck2020MLPelliptic} on MLP approximation algorithms contain proofs that the proposed MLP approximation algorithms overcome the curse of dimensionality for various types of nonlinear PDEs and thereby established, for the first time, that semilinear PDEs can actually be approximated without the curse of dimensionality.
However, none of these articles contain numerical simulations.
It is the subject of this article to generalize the MLP approximation algorithms in \cite{OvercomingPaper,Beck2020MLPAC,VW2019defaultRisk} to systems of PDEs and to present numerical simulations for several example PDEs in up to 1000 dimensions.
More precisely, 
	in \cref{sect:algorithm} we specify the generalized MLP approximation scheme which we propose in this article (see \eqref{setting:MLP} in \cref{sect:algorithm} below) and
	in \cref{sect:num_examples} we apply this numerical approximation scheme to four different kinds of semilinear PDEs. 
We consider 
	Allen-Cahn PDEs in \cref{subsect:AllenCahn},
	Sine-Gordon type PDEs in \cref{subsect:SineGordon},
	systems of coupled semilinear heat PDEs in \cref{subsect:PDESystem}, and
	semilinear Black-Scholes PDEs in  \cref{subsect:BlackScholes}.
In the case of each of the above mentioned example PDEs we approximately compute the relative $L^2$-error of the proposed MLP approximation algorithm 
(see 
	\cref{AllenCahn-MLP} and \cref{fig:AllenCahn} in \cref{subsect:AllenCahn},
	\cref{SineGordon-MLP} and \cref{fig:SineGordon} in \cref{subsect:SineGordon},
	\cref{PDESystem-MLP} and \cref{fig:PDESystem} in \cref{subsect:PDESystem}, and
	\cref{BlackScholes-MLP} and \cref{fig:BlackScholes} in \cref{subsect:BlackScholes}).
In our approximate computations of the relative $L^2$-errors the unknown exact solutions of the PDEs have been approximated by means of the deep learning based approximation method in Beck et al.\ \cite{beck2019deep}, the so-called deep splitting (DS) method 
(see the
	\nth{4} and \nth{5} columns in \cref{AllenCahn-MLP,SineGordon-MLP,PDESystem-MLP,BlackScholes-MLP} and 
	\cref{fig:AllenCahn_DS,fig:SineGordon_DS,fig:PDESystem_DS,fig:BlackScholes_DS})
and by means of the MLP approximation algorithm itself 
(see the
	\nth{4} and \nth{6} columns in \cref{AllenCahn-MLP,SineGordon-MLP,PDESystem-MLP,BlackScholes-MLP} and 
	\cref{fig:AllenCahn_MLP,fig:SineGordon_MLP,fig:PDESystem_MLP,fig:BlackScholes_MLP}).
In \cref{sect:code} we present the \Cplusplus \, source code employed to perform the numerical simulations presented in \cref{sect:num_examples}.

%%%%%%%%%%%%%%%%%%%%%%%%%%%%%%%%%%%%%%%%%%%
%%%%%%%%%%%%%%%%%%%%%%%%%%%%%%%%%%%%%%%%%%%
\section{Description of the approximation algorithm} 
\label{sect:algorithm}
%%%%%%%%%%%%%%%%%%%%%%%%%%%%%%%%%%%%%%%%%%%
%%%%%%%%%%%%%%%%%%%%%%%%%%%%%%%%%%%%%%%%%%%

In this section we introduce the generalized MLP approximation scheme which we consider in this article (see \eqref{setting:MLP} in \cref{nogradsetting} below).

\begin{setting}\label{nogradsetting}
Let 
$d,k\in \N$, 
$c,T\in (0,\infty)$, 
$\Theta = \cup_{n\in \N} \Z^n$, 
$f=(f_{1},\dots,f_{k})\in C( \R^d \times \R^k,\R^k)$, 
$g\in C(\R^d,\R^k)$, 
let $\norm{\cdot} \colon (\cup_{q\in\N}\R^{q})\to [0,\infty)$ be the standard norm, 
let $\phi_r\colon \R^k \rightarrow \R^k$, 
	$r\in [0,\infty]$, 
	satisfy for all $r\in [0,\infty]$, $y = (y_1,\dots,y_k)\in \R^k$ that
\begin{equation}\label{algo:cut-off}
	\phi_r(y) = \left(\min\{r,\max\{-r,y_1\}\}, \dots, \min\{r,\max\{-r,y_k\}\} \right),
\end{equation}
let $(\Omega,\mathcal{F}, \P)$ be a probability space,
let $\mathcal{R}^\theta\colon \Omega \rightarrow [0,1],$ $\theta \in \Theta$, be independent $\mathcal{U}_{[0,1]}$-distributed random variables, 
let $W^\theta\colon [0,T]\times \Omega \rightarrow \R^d$, $\theta \in \Theta$, be independent standard Brownian motions, 
assume that $(\mathcal{R}^\theta)_{\theta \in \Theta}$ and $(W^\theta)_{\theta \in \Theta}$ are independent, 
let  $R^\theta\colon [0,T]\times \Omega \rightarrow [0,T]$, $\theta \in \Theta$, 
	satisfy for all $t\in [0,T]$, $\theta \in \Theta$ 
	that $R^\theta_{t} = t+ (T-t)\mathcal{R}^\theta$, 
let $\mu \colon \R^{d} \to \R^{d}$ 
	and $\sigma \colon \R^{d}\to \R^{d\times d}$ be globally Lipschitz continuous functions, 
for every 
	$x\in \R^{d}$, $\theta \in \Theta$, $t\in [0,T]$
let $(X^{x,\theta}_{t,s})_{s\in [t,T]} \colon [t,T] \times  \Omega \rightarrow \R^d$
	be a stochastic process with continuous sample paths 
	which satisfies that 
	for all $s\in [t,T]$ 
	it holds $\P$-a.s.\ that
\begin{equation}\label{setting:Xprocess}
	X^{x,\theta}_{t,s} 
	= x + \int_{t}^{s} \mu\big(X^{x,\theta}_{t,r}\big) \, dr + \int^{s}_{t} \sigma \big(X^{x,\theta}_{t,r}\big) \, dW^{\theta}_{r}, 
\end{equation}
let $V^\theta_{n,M,r}\colon [0,T]\times \R^d\times \Omega \rightarrow \R^k$, 
	$\theta \in \Theta$, 
	$n \in \Z$, 
	$M\in \N$, 
	$r\in [0,\infty]$, 
	satisfy for all 
	$\theta \in \Theta$, 
	$n, M \in \N$, 
	$r\in [0,\infty]$, 
	$t\in [0,T]$, 
	$x\in \R^d$ 
	that 
	$V^\theta_{-1,M,r}(t,x) = V^\theta_{0,M,r}(t,x) = 0$ and 
\begin{equation}
\label{setting:MLP}
\begin{split}
	V^\theta_{n,M,r}(t,x) 
		&= \frac{1}{M^n} \Bigg[ \,\sum_{m=1}^{M^n} g\big(X^{x,(\theta,0,-m)}_{t,T}\big)  \Bigg]\\
		&+ \sum_{l=0}^{n-1} \frac{(T-t)}{M^{n-l}}  
			\sum_{m=1}^{M^{n-l}} 
			 \bigg[ f \Big(X^{x,(\theta,l,m)}_{t,R_t^{(\theta,l,m)}},\phi_{r}\Big(V^{(\theta,l,m)}_{l,M,r}\big(R_t^{(\theta,l,m)},X^{x,(\theta,l,m)}_{t,R_t^{(\theta,l,m)}}\big)\Big)\Big) \\
		&\qquad- \mathbbm{1}_\N(l) \, f \Big(X^{x,(\theta,l,m)}_{t,R_t^{(\theta,l,m)}},\phi_{r}\Big(V^{(\theta,l,-m)}_{l-1,M,r}\big(R_t^{(\theta,l,m)},X^{x,(\theta,l,m)}_{t,R_t^{(\theta,l,m)}}\big)\Big)\Big)\bigg] ,
\end{split}
\end{equation}
and let $u = (u(t, (x_1, x_2, \ldots, x_d)))_{(t, (x_1, x_2, \ldots, x_d)) \in [0,T] \times \R^d} =(u_{1},\dots,u_{k})  \in C([0,T]\times \R^d, \R^{k})$ satisfy for all 
	$t\in [0,T)$, 
	$x \in\R^d$, $i\in \{1,2,\dots,k\}$ 
	that 
	$u\vert_{[0,T)\times \R^d}\in C^{1,2}([0,T)\times \R^d,\R^{k})$, 
	$\norm{u(t,x)}\leq c(1+\norm{x}^{c})$, 
	$u(T,x) = g(x)$, 
	and
	
\begin{equation}
\label{setting:PDE}
\begin{split} 
	&(\tfrac{\partial}{\partial t}u_{i})(t,x) + \tfrac 12 \Trace\big(\sigma(x)[\sigma(x)]^{*}(\Hess_{x} u_{i})(t,x)\big) \\
	&+\bigg[\smallsum\limits_{j=1}^{d}\mu_{j}(x) (\tfrac{\del}{\del x_{j}}u_{i})(t,x)\bigg] + f_{i}(x,u(t,x))
=
	0.
\end{split}
\end{equation}
\end{setting}

%%%%%%%%%%%%%%%%%%%%%%%%%%%%%%%%%%%%%%%%%%%
%%%%%%%%%%%%%%%%%%%%%%%%%%%%%%%%%%%%%%%%%%%
\section{Numerical examples}
\label{sect:num_examples}
%%%%%%%%%%%%%%%%%%%%%%%%%%%%%%%%%%%%%%%%%%%
%%%%%%%%%%%%%%%%%%%%%%%%%%%%%%%%%%%%%%%%%%%

In this section we apply the MLP approximation algorithm in \eqref{setting:MLP} in \cref{nogradsetting} above to four different kinds of semilinear PDEs. 
We consider 
	Allen-Cahn PDEs in \cref{subsect:AllenCahn},
	Sine-Gordon type PDEs in \cref{subsect:SineGordon},
	systems of coupled semilinear heat PDEs in \cref{subsect:PDESystem}, and
	semilinear Black-Scholes PDEs in  \cref{subsect:BlackScholes}.

%%%%%%%%%%%%%%%%%%%%%%%%%%%%%%%%%%%%%%%%%%%
\subsection{Allen-Cahn partial differential equations (PDEs)}
\label{subsect:AllenCahn}
%%%%%%%%%%%%%%%%%%%%%%%%%%%%%%%%%%%%%%%%%%%

In this subsection we apply the MLP approximation algorithm in \eqref{setting:MLP} in \cref{nogradsetting} above to the Allen-Cahn PDE in \eqref{AllenCahn:eq1} below (cf., e.g., Bartels \cite[Chapter 6]{Bartels2015} and Feng \& Prohl \cite{Feng2003}).

\begin{figure} 
\centering 
\begin{subfigure}{0.49\textwidth}
    \includegraphics[width=\textwidth]{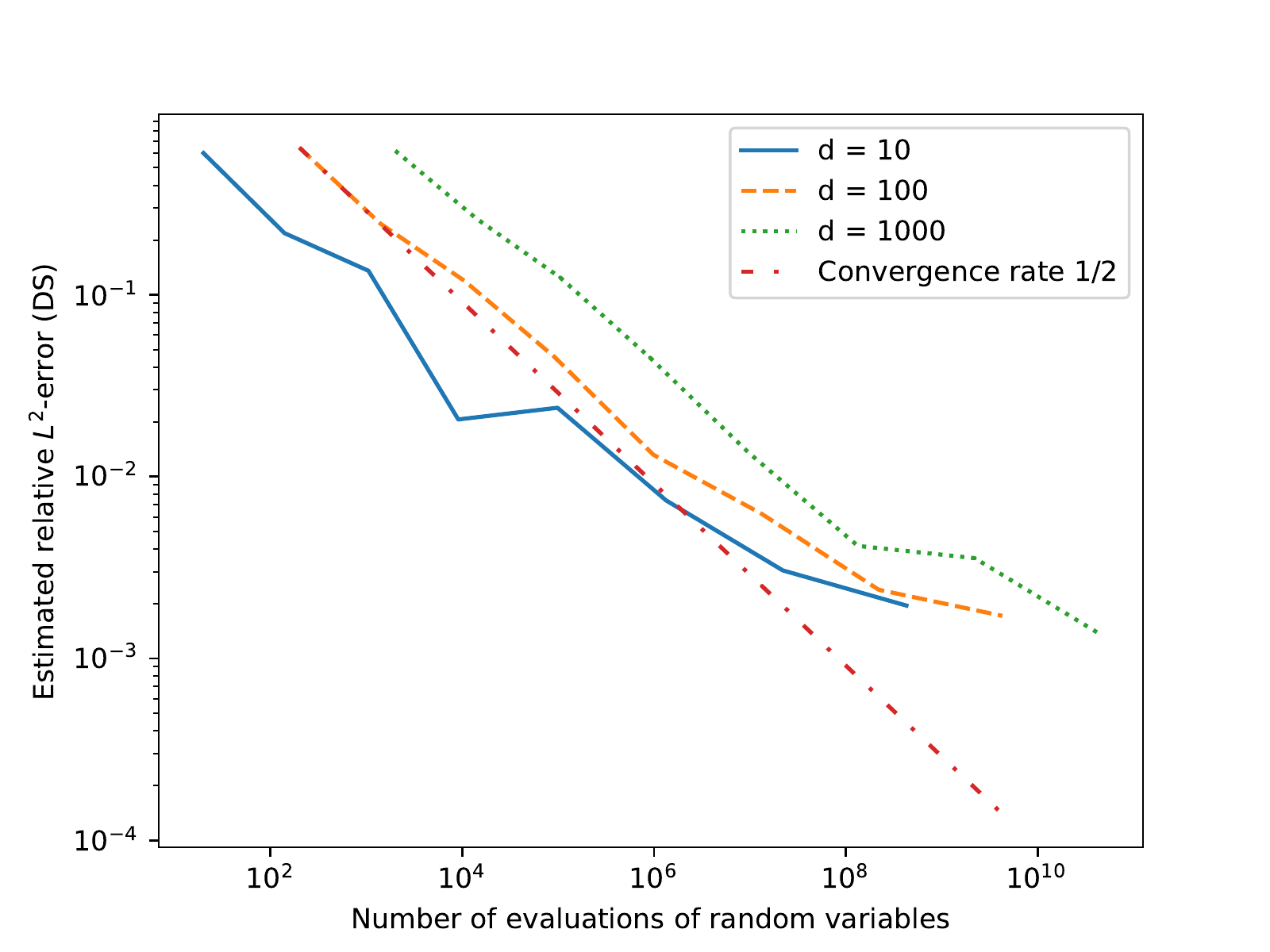} 
    \caption{Reference solutions computed by DS} 
    \label{fig:AllenCahn_DS} 
\end{subfigure} 
\begin{subfigure}{0.49\textwidth} 
    \includegraphics[width=\textwidth]{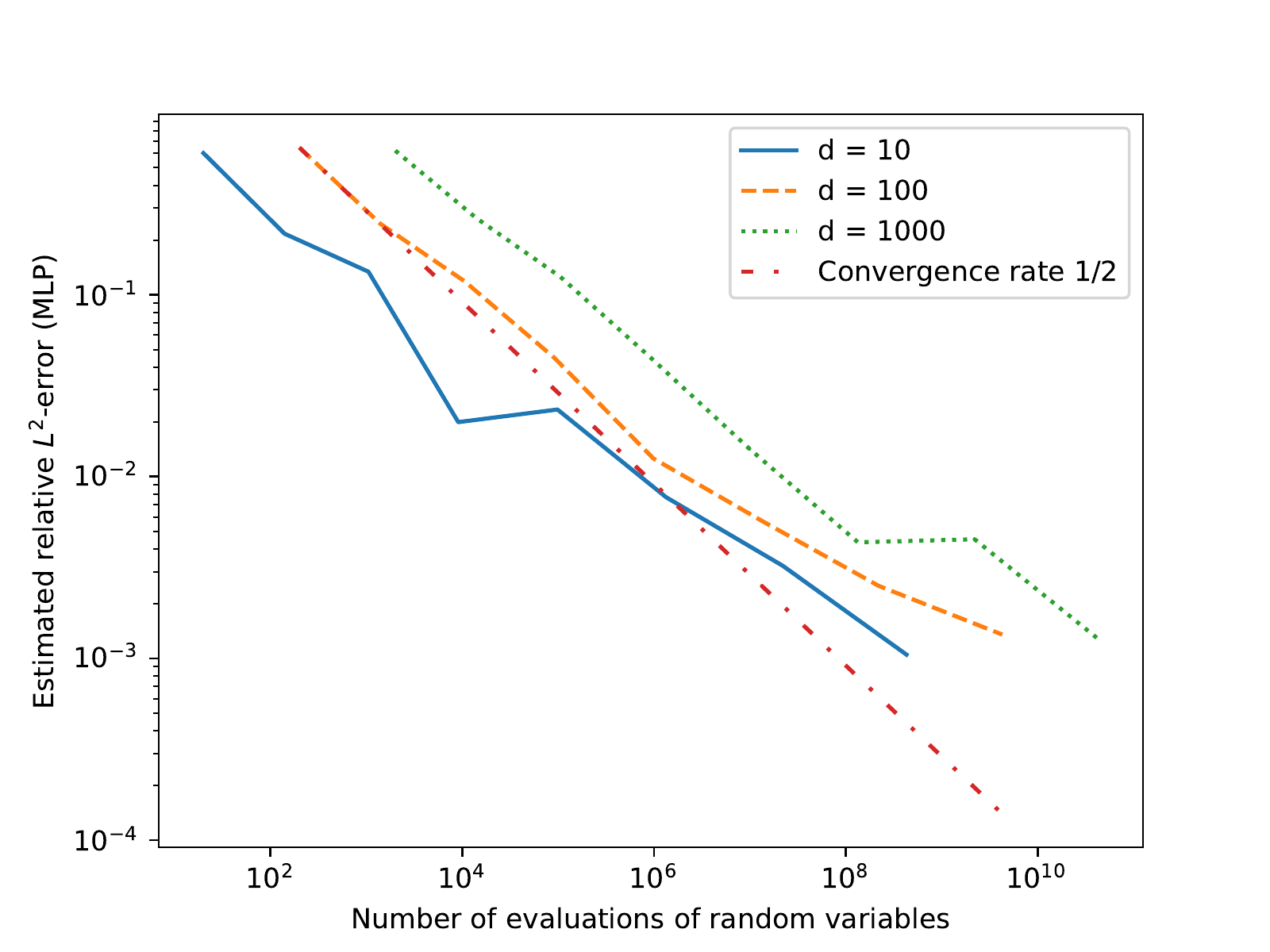} 
    \caption{Reference solutions computed by MLP}
    \label{fig:AllenCahn_MLP}
\end{subfigure}
\caption{Approximative plot of the relative $L^2$-error of the MLP approximation algorithm in \eqref{setting:MLP} against the computational effort of the algorithm in the case of the Allen-Cahn PDE in \eqref{AllenCahn:eq1}.} 
\label{fig:AllenCahn} 
\end{figure}

Assume  \cref{nogradsetting} 
and assume for all 
$t\in [0,T]$, 
$x,v\in \R^d$, 
$y\in \R$ 
that 
$k=1$, 
$T=1$, 
$f(x,y)= y-y^3$, 
$g(x) = (2+\tfrac 25 \norm{x}^{2})^{-1}$, 
$\mu(x)=0$, 
and 
$\sigma(x)v = \sqrt{2}v$. 
Note that this, (\ref{setting:Xprocess}), and (\ref{setting:PDE}) 
assure that 
for all 
$x\in \R^d$, 
$\theta \in \Theta$, 
$t\in [0,T]$, $s\in[t,T]$ 
it holds that
 $\P\big(X^{x,\theta}_{t,s} = x+\sqrt{2}(W^{\theta}_{s}-W^{\theta}_{t})\big)=1$ 
 and
 \begin{equation}
 \label{AllenCahn:eq1}
(\tfrac{\partial}{\partial t}u)(t,x) + (\Delta_{x}u)(t,x) + u(t,x)-(u(t,x))^{3}= 0.
\end{equation}
Observe that 
for all $x\in\R^{d}$, $y\in\R$
it holds that 
$yf(x,y)= y^{2}-y^{4} \leq 1+ y^{2}$. 
Combining this and (\ref{setting:PDE}) with Beck et al.\ \cite[Corollary 2.4]{Beck2020MLPAC} 
ensures that
for all $t\in [0,T]$ 
 it holds that
% \todo{Explain why we chose $r = 4$.}
\begin{equation}
\sup_{x\in\R^{d}} |u(t,x)| 
\leq e^{T-t}\big[1+\sup_{x\in\R^{d}}|u(T,x)|^{2}\big]^{\!\nicefrac 12} 
%= e^{T-t}\big[1+\sup_{x\in\R^{d}}|g(x)|^{2}\big]^{\!\nicefrac 12} 
\leq e \big[1 + \tfrac{1}{4}\big]^{\!\nicefrac 12} 
= \frac{\sqrt 5 e } 2 \leq 4.
\end{equation}
%
%
%%%%%%%%%%%%%%%%%%%%%%%%%%%%%%%%%%%%%%%%% 
%START 
%\input{01_Figure_descriptions/AllenCahn_description.tex}
%%%%%%%%%%%%%%%%%%%%%%%%%%%%%%%%%%%%%%%%% 
In \cref{AllenCahn-MLP} we approximately present for 
    $d \in \{10, 100, 1000 \}$, 
    $n \in \{1, 2, \ldots, 8 \}$ 
one random realization of $V^\theta_{n,n,4}(0,0)$ 
    (\nth {3} column in \cref{AllenCahn-MLP}), 
the relative $L^2$-error $\frac{( \E[ |V^{(0)}_{n,n,4}(0,0)-u(0,0)|^{2}])^{1/2}}{u(0,0)}$
    (\nth {5} and \nth {6} column in \cref{AllenCahn-MLP}),
the number of evaluations of one-dimensional random variables used to calculate one random realization of $V^\theta_{n,n,4}(0,0)$
    (\nth {7} column in \cref{AllenCahn-MLP}),
and
the runtime to calculate one random realization of $V^\theta_{n,n,4}(0,0)$
    (\nth {8} column in \cref{AllenCahn-MLP}).
In \cref{fig:AllenCahn} we approximately plot for
    $d \in \{10, 100, 1000 \}$,
     $n \in \{1, 2, \ldots, 8 \}$ 
the relative $L^2$-error $\frac{( \E[ |V^{(0)}_{n,n,4}(0,0)-u(0,0)|^{2}])^{1/2}}{u(0,0)}$ 
    (\nth {5} and \nth {6} column in \cref{AllenCahn-MLP})
against the number of evaluations of one-dimensional random variables used to calculate one random realization of $V^\theta_{n,n,4}(0,0)$
    (\nth {7} column in \cref{AllenCahn-MLP}).
The results in \cref{AllenCahn-MLP} and \cref{fig:AllenCahn} have been computed by means of \Cplusplus~code~\ref{code} in Section~\ref{sect:code} below.
For every 
    $n \in \{1, 2, \ldots, 8 \}$ 
for our approximative computations of the relative $L^2$-error $\frac{( \E[ |V^{(0)}_{n,n,4}(0,0)-u(0,0)|^{2}])^{1/2}}{u(0,0)}$ 
    (\nth {5} and \nth {6} column in \cref{AllenCahn-MLP})
the value $u(0,0)$ of the unknown exact solution 
in the relative $L^2$-error
has been approximated by means of an average of 5 independent runs of the deep splitting approximation method in Beck et al.\ \cite{beck2019deep}
    (\nth {5} column in \cref{AllenCahn-MLP})
and by means of an average of 5 independent evaluations of $V^{(0)}_{8,8,4}(0,0)$ 
    (\nth {6} column in \cref{AllenCahn-MLP}),
respectively,
and the expectation 
in the relative $L^2$-error 
has been approximated by means of Monte Carlo approximations involving 5 independent runs.

%%%%%%%%%%%%%%%%%%%%%%%%%%%%%%%%%%%%%%%%% 
%%END 
%\input{01_Figure_descriptions/AllenCahn_description.tex}
%%%%%%%%%%%%%%%%%%%%%%%%%%%%%%%%%%%%%%%%% 

%%%%%%%%%%%%%%%%%%%%%%%%%%%%%%%%%%%%%%%%%%%
%Numerical results
%%%%%%%%%%%%%%%%%%%%%%%%%%%%%%%%%%%%%%%%%%%

%%%%%%%%%%%%%%%%%%%%%%%%%%%%%%%%%%%%%%%%% 
%START 
%\input{../2_CSV_data/07_Simulations_200510_AC/1_Converted_files/Tabular/AllenCahn-MLP_tabular}
%%%%%%%%%%%%%%%%%%%%%%%%%%%%%%%%%%%%%%%%% 
\begin{table} 
\begin{center} 
\begin{tabular}{|c|c|c|c|c|c|c|c|}
\hline 
\thead{d} &\thead{n} &\thead{Result \\ of \\ MLP \\ algo- \\ rithm} &\thead{Refe- \\ rence \\ solu- \\ tions} &\thead{Esti- \\ mated \\ relative \\ $L^2$-error \\ (DS)} &\thead{Esti- \\ mated \\ relative \\ $L^2$-error \\ (MLP)} &\thead{Evaluations \\ of \\ random \\ variables} &\thead{Run- \\ time \\ in \\ sec- \\ onds} \\ 
\hline 
10 &1 &0.09925 &  &0.601730 &0.600980 &20 &0.00016 \\ 
10 &2 &0.22436 &DS: &0.218321 &0.216805 &140 &0.00012 \\ 
10 &3 &0.24525 &0.29614 &0.135595 &0.133916 &1050 &0.00059 \\ 
10 &4 &0.28409 &  &0.020630 &0.019931 &9080 &0.00323 \\ 
10 &5 &0.30594 &  &0.023919 &0.023377 &98300 &0.03291 \\ 
10 &6 &0.29642 &MLP: &0.007396 &0.007699 &1334340 &0.45886 \\ 
10 &7 &0.29662 &0.29555 &0.003047 &0.003228 &22032010 &4.16987 \\ 
10 &8 &0.29555 &  &0.001953 &0.001049 &428332080 &105.098 \\ 
\hline 
100 &1 &0.01350 &  &0.645486 &0.645186 &200 &0.00006 \\ 
100 &2 &0.02433 &DS: &0.246201 &0.245570 &1400 &0.00012 \\ 
100 &3 &0.03433 &0.03376 &0.119473 &0.118891 &10500 &0.00087 \\ 
100 &4 &0.03345 &  &0.045691 &0.045087 &90800 &0.00857 \\ 
100 &5 &0.03332 &  &0.013188 &0.012581 &983000 &0.09031 \\ 
100 &6 &0.03346 &MLP: &0.006225 &0.005701 &13343400 &1.29344 \\ 
100 &7 &0.03375 &0.03373 &0.002386 &0.002504 &220320100 &25.3006 \\ 
100 &8 &0.03373 &  &0.001714 &0.001351 &4283320800 &827.336 \\ 
\hline 
1000 &1 &0.00128 &  &0.620073 &0.620698 &2000 &0.00011 \\ 
1000 &2 &0.00253 &DS: &0.262782 &0.263995 &14000 &0.00063 \\ 
1000 &3 &0.00298 &0.00339 &0.124379 &0.125679 &105000 &0.00446 \\ 
1000 &4 &0.00366 &  &0.045096 &0.045392 &908000 &0.04711 \\ 
1000 &5 &0.00338 &  &0.013440 &0.014241 &9830000 &0.56212 \\ 
1000 &6 &0.00340 &MLP: &0.004158 &0.004351 &133434000 &7.77024 \\ 
1000 &7 &0.00340 &0.00340 &0.003561 &0.004518 &2203201000 &209.154 \\ 
1000 &8 &0.00340 &  &0.001379 &0.001278 &42833208000 &7786.58 \\ 
\hline 
\end{tabular} 
\end{center} 
\caption{ \label{AllenCahn-MLP} Numerical simulations for the MLP approximation algorithm in \eqref{setting:MLP} in the case of the Allen-Cahn PDE in \eqref{AllenCahn:eq1}} 
\end{table} 

%%%%%%%%%%%%%%%%%%%%%%%%%%%%%%%%%%%%%%%%% 
%%END 
%\input{../2_CSV_data/07_Simulations_200510_AC/1_Converted_files/Tabular/AllenCahn-MLP_tabular}
%%%%%%%%%%%%%%%%%%%%%%%%%%%%%%%%%%%%%%%%% 

%%%%%%%%%%%%%%%%%%%%%%%%%%%%%%%%%%%%%%%%% 
%START 
%\input{../2_CSV_data/07_Simulations_200510_AC/1_Converted_files/Plots/AllenCahn_plotfigure}
%%%%%%%%%%%%%%%%%%%%%%%%%%%%%%%%%%%%%%%%% 

%%%%%%%%%%%%%%%%%%%%%%%%%%%%%%%%%%%%%%%%% 
%%END 
%\input{../2_CSV_data/07_Simulations_200510_AC/1_Converted_files/Plots/AllenCahn_plotfigure}
%%%%%%%%%%%%%%%%%%%%%%%%%%%%%%%%%%%%%%%%% 

%%%%%%%%%%%%%%%%%%%%%%%%%%%%%%%%%%%%%%%%%%%
\subsection{Sine-Gordon type PDEs}
\label{subsect:SineGordon}
%%%%%%%%%%%%%%%%%%%%%%%%%%%%%%%%%%%%%%%%%%%

In this subsection we apply the MLP approximation algorithm in \eqref{setting:MLP} in \cref{nogradsetting} above to the Sine-Gordon-type PDE in \eqref{SineGordon:eq1} below (cf., e.g., Hairer \& Hao \cite{shen2016sinegordon}, Barone \cite{barone1971theory}, and Coleman \cite{coleman1975gordon}).

\begin{figure} 
\centering 
\begin{subfigure}{0.49\textwidth}
    \includegraphics[width=\textwidth]{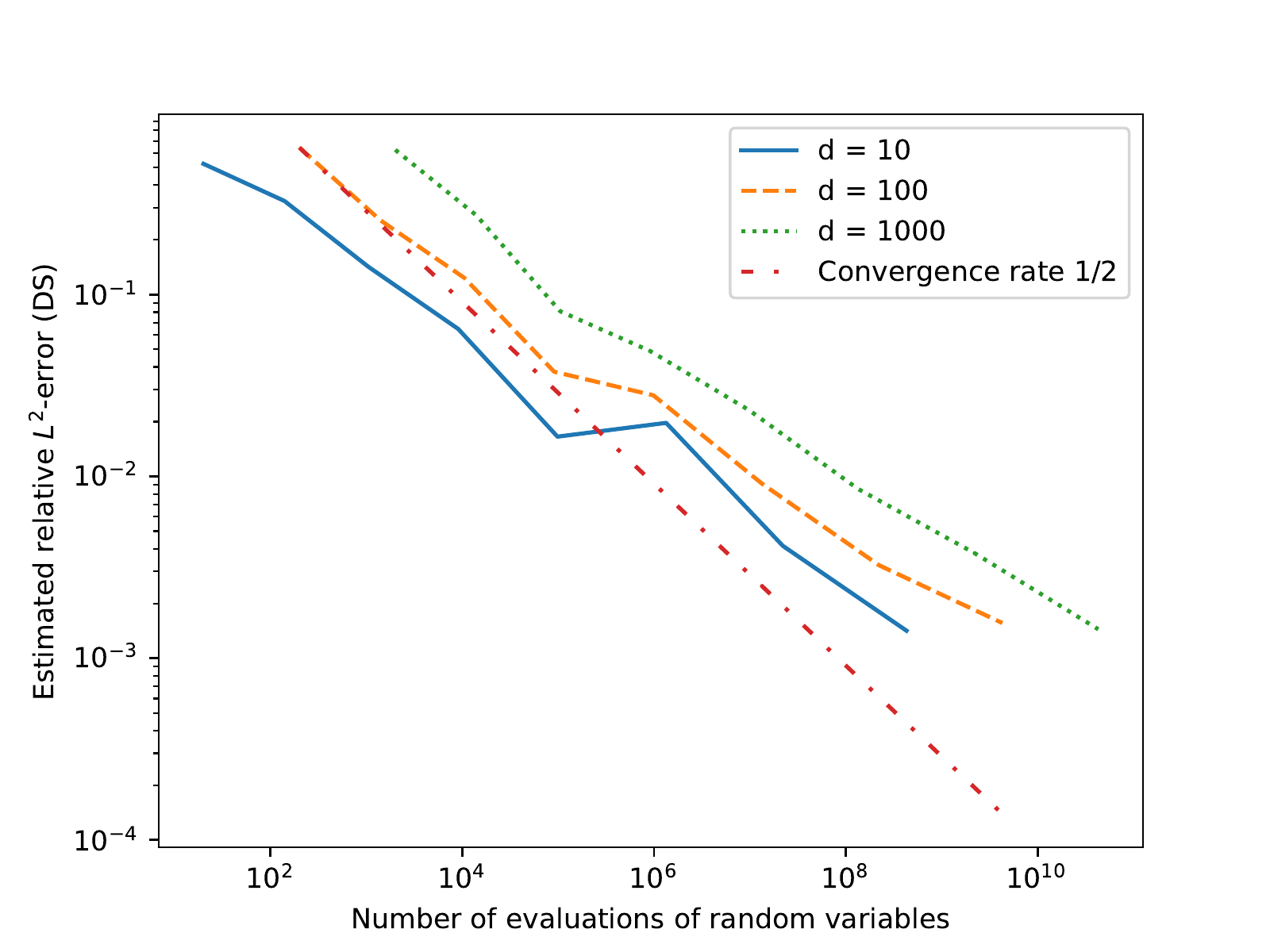} 
    \caption{Reference solutions computed by DS} 
    \label{fig:SineGordon_DS} 
\end{subfigure} 
\begin{subfigure}{0.49\textwidth} 
    \includegraphics[width=\textwidth]{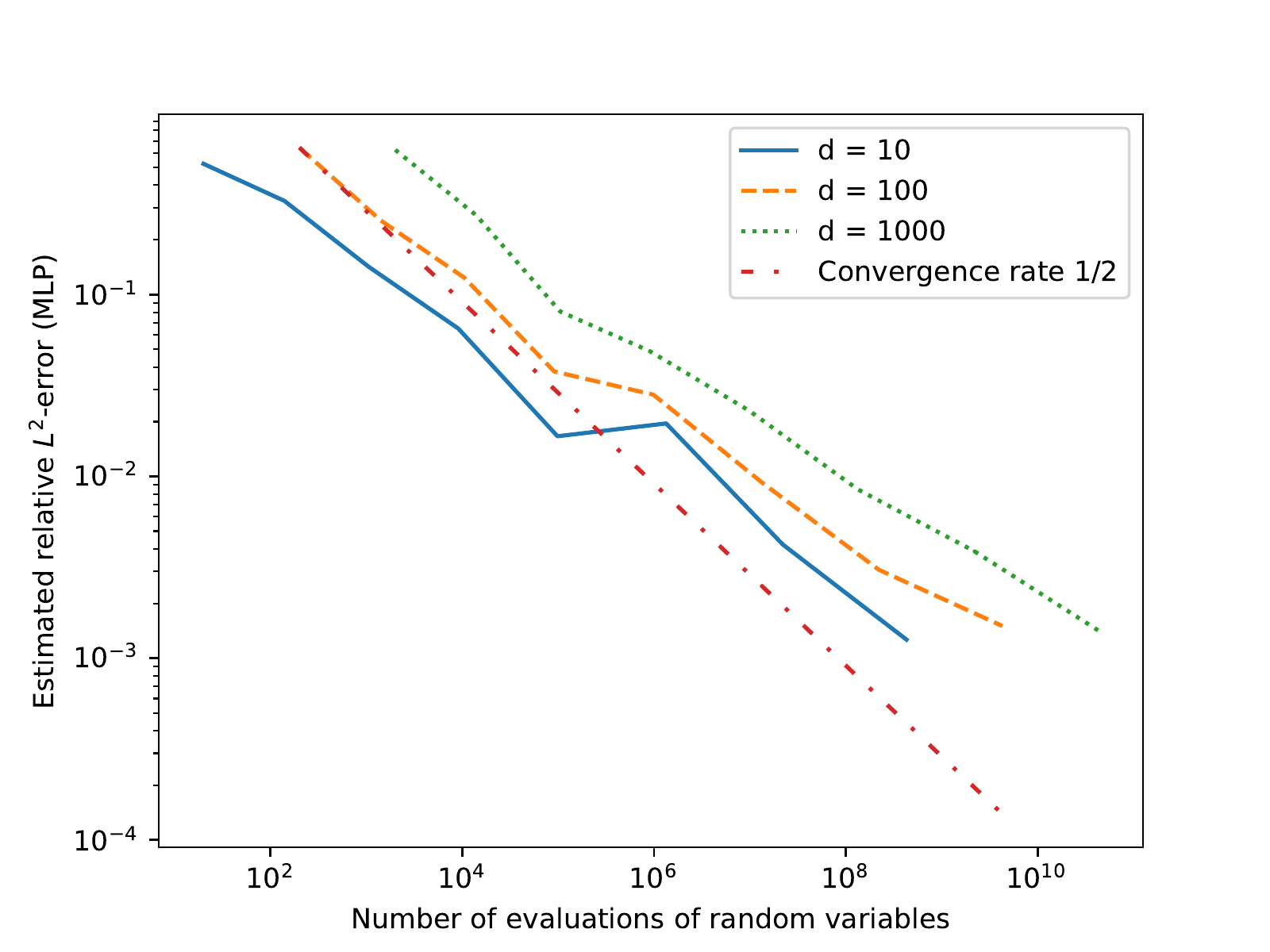} 
    \caption{Reference solutions computed by MLP}
    \label{fig:SineGordon_MLP}
\end{subfigure}
\caption{Approximative plot of the relative $L^2$-error of the MLP approximation algorithm in \eqref{setting:MLP} against the computational effort of the algorithm in the case of the Sine-Gordon-type PDE in \eqref{SineGordon:eq1}.} 
\label{fig:SineGordon} 
\end{figure}

Assume \cref{nogradsetting} 
and 
assume for all 
$t\in [0,T]$, 
$x,v\in \R^d$, 
$y\in \R$ 
that 
%$d = 100$, 
$k=1$, 
$T=1$,
$f(x,y)= \sin (y)$, 
$g(x) = (2+\tfrac 25 \norm{x}^{2})^{-1}$, 
$\mu(x)=0$, 
and 
$\sigma(x)v = \sqrt{2}v$. 
This, \eqref{setting:Xprocess}, and \eqref{setting:PDE}
ensure that for all 
$x\in \R^d$, 
$\theta \in \Theta$, 
$t\in [0,T]$, 
$s\in[t,T]$ 
it holds that 
$\P\big(X^{x,\theta}_{t,s} = x+\sqrt{2}(W^{\theta}_{s}-W^{\theta}_{t})\big)=1$ 
and
 \begin{equation}
 \label{SineGordon:eq1}
(\tfrac{\partial}{\partial t}u)(t,x) + (\Delta_{x}u)(t,x) +  \sin(u(t,x)) = 0.
\end{equation}
%
%
%%%%%%%%%%%%%%%%%%%%%%%%%%%%%%%%%%%%%%%%% 
%START 
%\input{01_Figure_descriptions/SineGordon_description.tex}
%%%%%%%%%%%%%%%%%%%%%%%%%%%%%%%%%%%%%%%%% 
In \cref{SineGordon-MLP} we approximately present for 
    $d \in \{10, 100, 1000 \}$, 
    $n \in \{1, 2, \ldots, 8 \}$ 
one random realization of $V^\theta_{n,n,\infty}(0,0)$ 
    (\nth {3} column in \cref{SineGordon-MLP}), 
the relative $L^2$-error $\frac{( \E[ |V^{(0)}_{n,n,\infty}(0,0)-u(0,0)|^{2}])^{1/2}}{u(0,0)}$
    (\nth {5} and \nth {6} column in \cref{SineGordon-MLP}),
the number of evaluations of one-dimensional random variables used to calculate one random realization of $V^\theta_{n,n,\infty}(0,0)$
    (\nth {7} column in \cref{SineGordon-MLP}),
and
the runtime to calculate one random realization of $V^\theta_{n,n,\infty}(0,0)$
    (\nth {8} column in \cref{SineGordon-MLP}).
In \cref{fig:SineGordon} we approximately plot for
    $d \in \{10, 100, 1000 \}$,
     $n \in \{1, 2, \ldots, 8 \}$ 
the relative $L^2$-error $\frac{( \E[ |V^{(0)}_{n,n,\infty}(0,0)-u(0,0)|^{2}])^{1/2}}{u(0,0)}$ 
    (\nth {5} and \nth {6} column in \cref{SineGordon-MLP})
against the number of evaluations of one-dimensional random variables used to calculate one random realization of $V^\theta_{n,n,\infty}(0,0)$
    (\nth {7} column in \cref{SineGordon-MLP}).
The results in \cref{SineGordon-MLP} and \cref{fig:SineGordon} have been computed by means of \Cplusplus~code~\ref{code} in Section~\ref{sect:code} below.
For every 
    $n \in \{1, 2, \ldots, 8 \}$ 
for our approximative computations of the relative $L^2$-error $\frac{( \E[ |V^{(0)}_{n,n,\infty}(0,0)-u(0,0)|^{2}])^{1/2}}{u(0,0)}$ 
    (\nth {5} and \nth {6} column in \cref{SineGordon-MLP})
the value $u(0,0)$ of the unknown exact solution 
in the relative $L^2$-error
has been approximated by means of an average of 5 independent runs of the deep splitting approximation method in Beck et al.\ \cite{beck2019deep}
    (\nth {5} column in \cref{SineGordon-MLP})
and by means of an average of 5 independent evaluations of $V^{(0)}_{8,8,\infty}(0,0)$ 
    (\nth {6} column in \cref{SineGordon-MLP}),
respectively,
and the expectation 
in the relative $L^2$-error 
has been approximated by means of Monte Carlo approximations involving 5 independent runs.

%%%%%%%%%%%%%%%%%%%%%%%%%%%%%%%%%%%%%%%%% 
%%END 
%\input{01_Figure_descriptions/SineGordon_description.tex}
%%%%%%%%%%%%%%%%%%%%%%%%%%%%%%%%%%%%%%%%% 

%%%%%%%%%%%%%%%%%%%%%%%%%%%%%%%%%%%%%%%%%%%
%Numerical results
%%%%%%%%%%%%%%%%%%%%%%%%%%%%%%%%%%%%%%%%%%%

%%%%%%%%%%%%%%%%%%%%%%%%%%%%%%%%%%%%%%%%% 
%START 
%\input{../2_CSV_data/06_Simulations_200310/Merged_csv_files/Tabular/SineGordon-MLP_tabular}
%%%%%%%%%%%%%%%%%%%%%%%%%%%%%%%%%%%%%%%%% 
\begin{table} 
\begin{center} 
\begin{tabular}{|c|c|c|c|c|c|c|c|}
\hline 
\thead{d} &\thead{n} &\thead{Result \\ of \\ MLP \\ algo- \\ rithm} &\thead{Refe- \\ rence \\ solu- \\ tions} &\thead{Esti- \\ mated \\ relative \\ $L^2$-error \\ (DS)} &\thead{Esti- \\ mated \\ relative \\ $L^2$-error \\ (MLP)} &\thead{Evaluations \\ of \\ random \\ variables} &\thead{Run- \\ time \\ in \\ sec- \\ onds} \\ 
\hline 
10 &1 &0.16709 &  &0.523870 &0.524147 &20 &0.00005 \\ 
10 &2 &0.23704 &DS: &0.327063 &0.327472 &140 &0.00007 \\ 
10 &3 &0.28555 &0.30603 &0.142115 &0.142569 &1050 &0.00045 \\ 
10 &4 &0.28834 &  &0.064670 &0.065077 &9080 &0.00287 \\ 
10 &5 &0.31199 &  &0.016543 &0.016628 &98300 &0.03513 \\ 
10 &6 &0.30894 &MLP: &0.019695 &0.019559 &1334340 &0.57765 \\ 
10 &7 &0.30453 &0.30623 &0.004147 &0.004217 &22032010 &4.61993 \\ 
10 &8 &0.30580 &  &0.001417 &0.001266 &428332080 &103.820 \\ 
\hline 
100 &1 &0.01383 &  &0.643679 &0.643825 &200 &0.00008 \\ 
100 &2 &0.02576 &DS: &0.256185 &0.256491 &1400 &0.00014 \\ 
100 &3 &0.03558 &0.03375 &0.123576 &0.123843 &10500 &0.00096 \\ 
100 &4 &0.03496 &  &0.037638 &0.037774 &90800 &0.00732 \\ 
100 &5 &0.03330 &  &0.027919 &0.028147 &983000 &0.08787 \\ 
100 &6 &0.03396 &MLP: &0.009120 &0.009256 &13343400 &1.43123 \\ 
100 &7 &0.03398 &0.03376 &0.003250 &0.003068 &220320100 &24.8743 \\ 
100 &8 &0.03383 &  &0.001561 &0.001504 &4283320800 &823.935 \\ 
\hline 
1000 &1 &0.00130 &  &0.625283 &0.625188 &2000 &0.00013 \\ 
1000 &2 &0.00247 &DS: &0.272910 &0.272726 &14000 &0.00061 \\ 
1000 &3 &0.00321 &0.00339 &0.080925 &0.080727 &105000 &0.00465 \\ 
1000 &4 &0.00338 &  &0.048937 &0.048760 &908000 &0.04418 \\ 
1000 &5 &0.00335 &  &0.023106 &0.022971 &9830000 &0.52934 \\ 
1000 &6 &0.00339 &MLP: &0.008539 &0.008476 &133434000 &7.63818 \\ 
1000 &7 &0.00341 &0.00339 &0.003780 &0.003868 &2203201000 &206.351 \\ 
1000 &8 &0.00339 &  &0.001443 &0.001421 &42833208000 &7835.75 \\ 
\hline 
\end{tabular} 
\end{center} 
\caption{ \label{SineGordon-MLP} Numerical simulations for the MLP approximation algorithm in \eqref{setting:MLP} in the case of the Sine-Gordon-type PDE in \eqref{SineGordon:eq1}} 
\end{table} 

%%%%%%%%%%%%%%%%%%%%%%%%%%%%%%%%%%%%%%%%% 
%%END 
%\input{../2_CSV_data/06_Simulations_200310/Merged_csv_files/Tabular/SineGordon-MLP_tabular}
%%%%%%%%%%%%%%%%%%%%%%%%%%%%%%%%%%%%%%%%% 

%%%%%%%%%%%%%%%%%%%%%%%%%%%%%%%%%%%%%%%%% 
%START 
%\input{../2_CSV_data/06_Simulations_200310/Merged_csv_files/Plots/SineGordon_plotfigure}
%%%%%%%%%%%%%%%%%%%%%%%%%%%%%%%%%%%%%%%%% 

%%%%%%%%%%%%%%%%%%%%%%%%%%%%%%%%%%%%%%%%% 
%%END 
%\input{../2_CSV_data/06_Simulations_200310/Merged_csv_files/Plots/SineGordon_plotfigure}
%%%%%%%%%%%%%%%%%%%%%%%%%%%%%%%%%%%%%%%%% 

%%%%%%%%%%%%%%%%%%%%%%%%%%%%%%%%%%%%%%%%%%%
\subsection{System of semilinear heat PDEs}
\label{subsect:PDESystem}
%%%%%%%%%%%%%%%%%%%%%%%%%%%%%%%%%%%%%%%%%%%

In this subsection we apply the MLP approximation algorithm in \eqref{setting:MLP} in \cref{nogradsetting} above to the system of coupled semilinear heat PDEs in \eqref{PDESystem:eq1} below.

\begin{figure} 
\centering 
\begin{subfigure}{0.49\textwidth}
    \includegraphics[width=\textwidth]{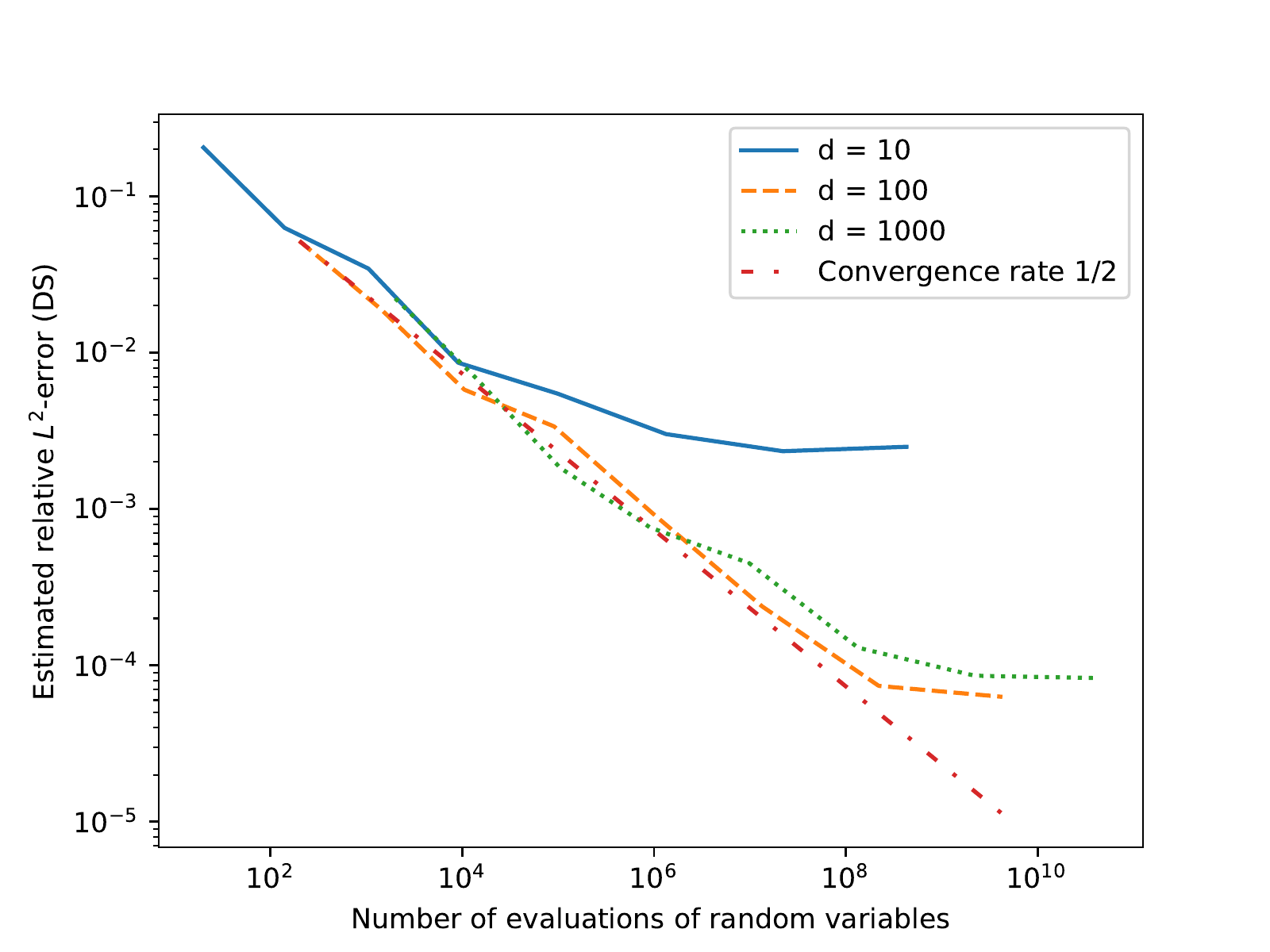} 
    \caption{Reference solutions computed by DS} 
    \label{fig:PDESystem_DS} 
\end{subfigure} 
\begin{subfigure}{0.49\textwidth} 
    \includegraphics[width=\textwidth]{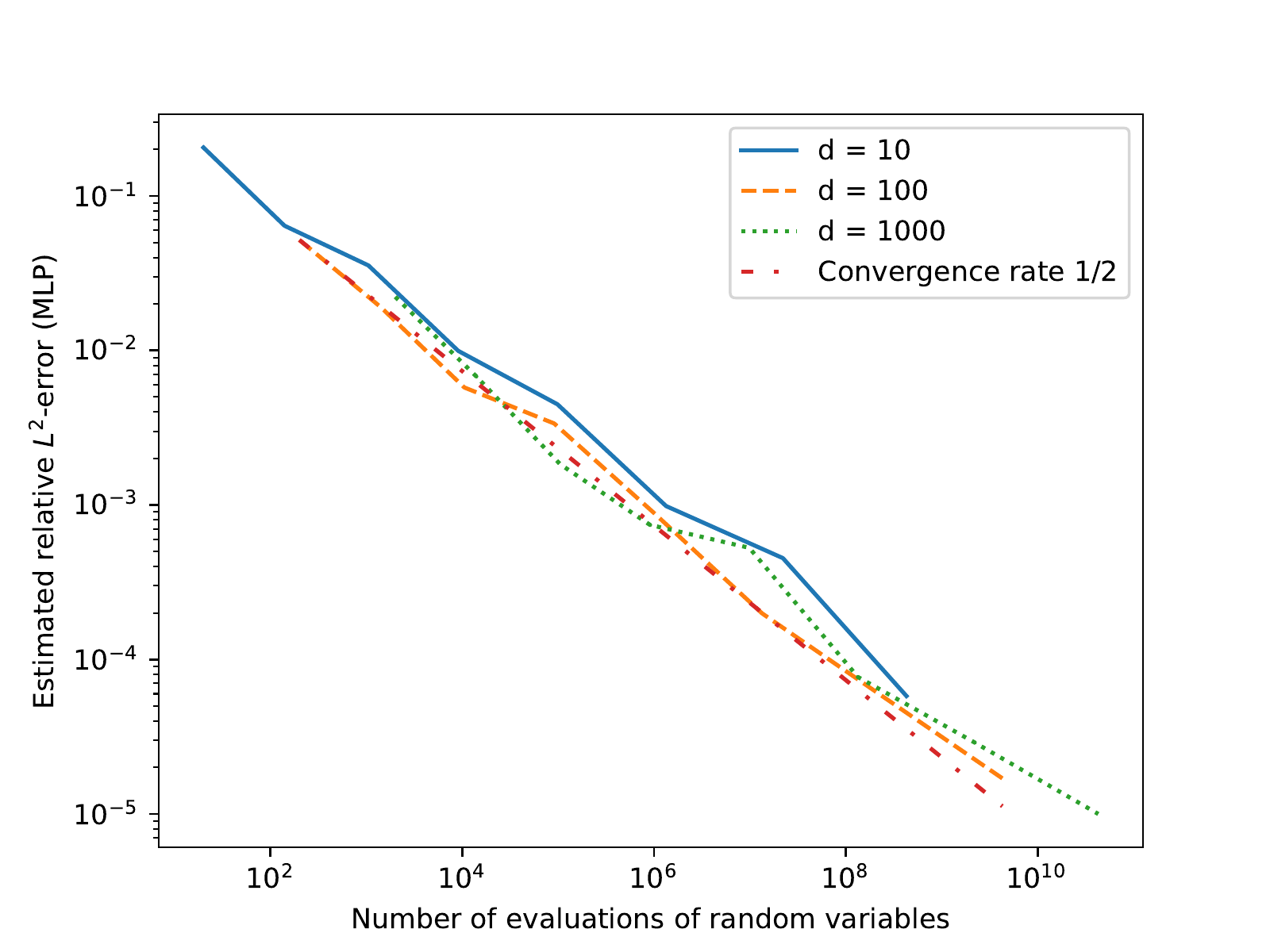} 
    \caption{Reference solutions computed by MLP}
    \label{fig:PDESystem_MLP}
\end{subfigure}
\caption{Approximative plot of the relative $L^2$-error of the MLP approximation algorithm in \eqref{setting:MLP} against the computational effort of the algorithm in the case of the system of coupled semilinear heat PDEs in \eqref{PDESystem:eq1}.} 
\label{fig:PDESystem} 
\end{figure}

Assume \cref{nogradsetting} 
and 
assume for all $t\in [0,T]$, $x,v\in \R^d$, $y=(y_{1},y_{2})\in \R^{2}$ 
that 
$k=2$, 
$T=1$, 
$f(x,y)= \big(\tfrac{y_{2}}{1+|y_{2}|^{2}},\tfrac {2 y_{1}}{3}\big)$, 
$g(x) = \big((2+\tfrac 25 \norm{x}^{2})^{-1}, \log(\tfrac 12 [1+\norm{x}^{2}])\big)$, 
$\mu(x)=0$, 
and 
$\sigma(x)v =\sqrt{2}v$. 
Observe that this, (\ref{setting:Xprocess}), and (\ref{setting:PDE}) implies that for all $x\in \R^d$, $\theta \in \Theta$, $t\in [0,T]$, $s\in[t,T]$ it holds that $\P\big(X^{x,\theta}_{t,s} = x+\sqrt{2}(W^{\theta}_{s}-W^{\theta}_{t})\big)=1$ and
\begin{equation}
\label{PDESystem:eq1}
(\tfrac{\partial}{\partial t}u)(t,x) +(\Delta_{x}u)(t,x)+ f\big(x,u(t,x)\big) = 0.
\end{equation}
%
%
%%%%%%%%%%%%%%%%%%%%%%%%%%%%%%%%%%%%%%%%% 
%START 
%\input{01_Figure_descriptions/PDESystem_description.tex}
%%%%%%%%%%%%%%%%%%%%%%%%%%%%%%%%%%%%%%%%% 
In \cref{PDESystem-MLP} we approximately present for 
    $d \in \{10, 100, 1000 \}$, 
    $n \in \{1, 2, \ldots, 8 \}$ 
one random realization of $V^\theta_{n,n,\infty}(0,0)$ 
    (\nth {3} column in \cref{PDESystem-MLP}), 
the relative $L^2$-error $\frac{( \E[ |V^{(0)}_{n,n,\infty}(0,0)-u(0,0)|^{2}])^{1/2}}{u(0,0)}$
    (\nth {5} and \nth {6} column in \cref{PDESystem-MLP}),
the number of evaluations of one-dimensional random variables used to calculate one random realization of $V^\theta_{n,n,\infty}(0,0)$
    (\nth {7} column in \cref{PDESystem-MLP}),
and
the runtime to calculate one random realization of $V^\theta_{n,n,\infty}(0,0)$
    (\nth {8} column in \cref{PDESystem-MLP}).
In \cref{fig:PDESystem} we approximately plot for
    $d \in \{10, 100, 1000 \}$,
     $n \in \{1, 2, \ldots, 8 \}$ 
the relative $L^2$-error $\frac{( \E[ |V^{(0)}_{n,n,\infty}(0,0)-u(0,0)|^{2}])^{1/2}}{u(0,0)}$ 
    (\nth {5} and \nth {6} column in \cref{PDESystem-MLP})
against the number of evaluations of one-dimensional random variables used to calculate one random realization of $V^\theta_{n,n,\infty}(0,0)$
    (\nth {7} column in \cref{PDESystem-MLP}).
The results in \cref{PDESystem-MLP} and \cref{fig:PDESystem} have been computed by means of \Cplusplus~code~\ref{code} in Section~\ref{sect:code} below.
For every 
    $n \in \{1, 2, \ldots, 8 \}$ 
for our approximative computations of the relative $L^2$-error $\frac{( \E[ |V^{(0)}_{n,n,\infty}(0,0)-u(0,0)|^{2}])^{1/2}}{u(0,0)}$ 
    (\nth {5} and \nth {6} column in \cref{PDESystem-MLP})
the value $u(0,0)$ of the unknown exact solution 
in the relative $L^2$-error
has been approximated by means of an average of 5 independent runs of the deep splitting approximation method in Beck et al.\ \cite{beck2019deep}
    (\nth {5} column in \cref{PDESystem-MLP})
and by means of an average of 5 independent evaluations of $V^{(0)}_{8,8,\infty}(0,0)$ 
    (\nth {6} column in \cref{PDESystem-MLP}),
respectively,
and the expectation 
in the relative $L^2$-error 
has been approximated by means of Monte Carlo approximations involving 5 independent runs.

%%%%%%%%%%%%%%%%%%%%%%%%%%%%%%%%%%%%%%%%% 
%%END 
%\input{01_Figure_descriptions/PDESystem_description.tex}
%%%%%%%%%%%%%%%%%%%%%%%%%%%%%%%%%%%%%%%%% 

%%%%%%%%%%%%%%%%%%%%%%%%%%%%%%%%%%%%%%%%%%%
%Numerical results
%%%%%%%%%%%%%%%%%%%%%%%%%%%%%%%%%%%%%%%%%%%

%%%%%%%%%%%%%%%%%%%%%%%%%%%%%%%%%%%%%%%%% 
%START 
%\input{../2_CSV_data/06_Simulations_200310/Merged_csv_files/Tabular/PDESystem-MLP_tabular}
%%%%%%%%%%%%%%%%%%%%%%%%%%%%%%%%%%%%%%%%% 
\begin{table} 
\begin{center} 
\begin{tabular}{|c|c|c|c|c|c|c|c|}
\hline 
\thead{d} &\thead{n} &\thead{Result \\ of MLP \\ algorithm} &\thead{Refe- \\ rence \\ solu- \\ tions} &\thead{Esti- \\ mated \\ relative \\ $L^2$-error \\ (DS)} &\thead{Esti- \\ mated \\ relative \\ $L^2$-error \\ (MLP)} &\thead{Evaluations \\ of \\ random \\ variables} &\thead{Run- \\ time \\ in \\ sec- \\ onds} \\ 
\hline 
10 &1 &(0.16714, 1.70085) &  &0.205522 &0.205894 &20 &0.00028 \\ 
10 &2 &(0.45134, 2.49508) &DS: &0.062991 &0.064206 &140 &0.00011 \\ 
10 &3 &(0.48243, 2.49320) &(0.47606, &0.034551 &0.035524 &1050 &0.00050 \\ 
10 &4 &(0.47085, 2.43225) &2.45101) &0.008637 &0.009941 &9080 &0.00306 \\ 
10 &5 &(0.47780, 2.44860) &  &0.005483 &0.004484 &98300 &0.03706 \\ 
10 &6 &(0.47539, 2.46088) &MLP: &0.003019 &0.000984 &1334340 &0.57155 \\ 
10 &7 &(0.47618, 2.45634) &(0.47621, &0.002347 &0.000453 &22032010 &4.54895 \\ 
10 &8 &(0.47627, 2.45705) &2.45726) &0.002504 &0.000058 &428332080 &105.486 \\ 
\hline 
100 &1 &(0.01481, 4.41159) &  &0.051931 &0.051942 &200 &0.00006 \\ 
100 &2 &(0.21898, 4.72326) &DS: &0.019148 &0.019169 &1400 &0.00012 \\ 
100 &3 &(0.22044, 4.71319) &(0.21892, &0.005808 &0.005762 &10500 &0.00081 \\ 
100 &4 &(0.21877, 4.67656) &4.67722) &0.003379 &0.003370 &90800 &0.00688 \\ 
100 &5 &(0.21992, 4.68292) &  &0.000931 &0.000895 &983000 &0.08848 \\ 
100 &6 &(0.21905, 4.67801) &MLP: &0.000239 &0.000199 &13343400 &1.42682 \\ 
100 &7 &(0.21898, 4.67743) &(0.21895, &0.000074 &0.000060 &220320100 &24.9459 \\ 
100 &8 &(0.21893, 4.67757) &4.67750) &0.000063 &0.000017 &4283320800 &828.779 \\ 
\hline 
1000 &1 &(0.00127, 6.88867) &  &0.022286 &0.022261 &2000 &0.00013 \\ 
1000 &2 &(0.14283, 6.87252) &DS: &0.006858 &0.006795 &14000 &0.00063 \\ 
1000 &3 &(0.14315, 6.94626) &(0.14273, &0.001836 &0.001834 &105000 &0.00487 \\ 
1000 &4 &(0.14244, 6.96112) &6.95587) &0.000763 &0.000743 &908000 &0.04441 \\ 
1000 &5 &(0.14278, 6.95906) &  &0.000451 &0.000528 &9830000 &0.51808 \\ 
1000 &6 &(0.14273, 6.95507) &MLP: &0.000130 &0.000077 &133434000 &7.68460 \\ 
1000 &7 &(0.14273, 6.95499) &(0.14274, &0.000086 &0.000029 &2203201000 &208.043 \\ 
1000 &8 &(0.14273, 6.95535) &6.95530) &0.000083 &0.000010 &42833208000 &7867.29 \\ 
\hline 
\end{tabular} 
\end{center} 
\caption{ \label{PDESystem-MLP} Numerical simulations for the MLP approximation algorithm in \eqref{setting:MLP} in the case of the system of coupled semilinear heat PDEs in \eqref{PDESystem:eq1}} 
\end{table} 

%%%%%%%%%%%%%%%%%%%%%%%%%%%%%%%%%%%%%%%%% 
%%END 
%\input{../2_CSV_data/06_Simulations_200310/Merged_csv_files/Tabular/PDESystem-MLP_tabular}
%%%%%%%%%%%%%%%%%%%%%%%%%%%%%%%%%%%%%%%%% 

%%%%%%%%%%%%%%%%%%%%%%%%%%%%%%%%%%%%%%%%% 
%START 
%\input{../2_CSV_data/06_Simulations_200310/Merged_csv_files/Plots/PDESystem_plotfigure}
%%%%%%%%%%%%%%%%%%%%%%%%%%%%%%%%%%%%%%%%% 

%%%%%%%%%%%%%%%%%%%%%%%%%%%%%%%%%%%%%%%%% 
%%END 
%\input{../2_CSV_data/06_Simulations_200310/Merged_csv_files/Plots/PDESystem_plotfigure}
%%%%%%%%%%%%%%%%%%%%%%%%%%%%%%%%%%%%%%%%% 

%%%%%%%%%%%%%%%%%%%%%%%%%%%%%%%%%%%%%%%%%%%
\subsection{Semilinear Black-Scholes PDEs}
\label{subsect:BlackScholes}
%%%%%%%%%%%%%%%%%%%%%%%%%%%%%%%%%%%%%%%%%%%

In this subsection we apply the MLP approximation algorithm in \eqref{setting:MLP} in \cref{nogradsetting} above to the semilinear Black-Scholes PDE in \eqref{BlackScholes:eq1} below  (cf. Black \& Scholes \cite{BlackScholes73}).

\begin{figure} 
\centering 
\begin{subfigure}{0.49\textwidth}
    \includegraphics[width=\textwidth]{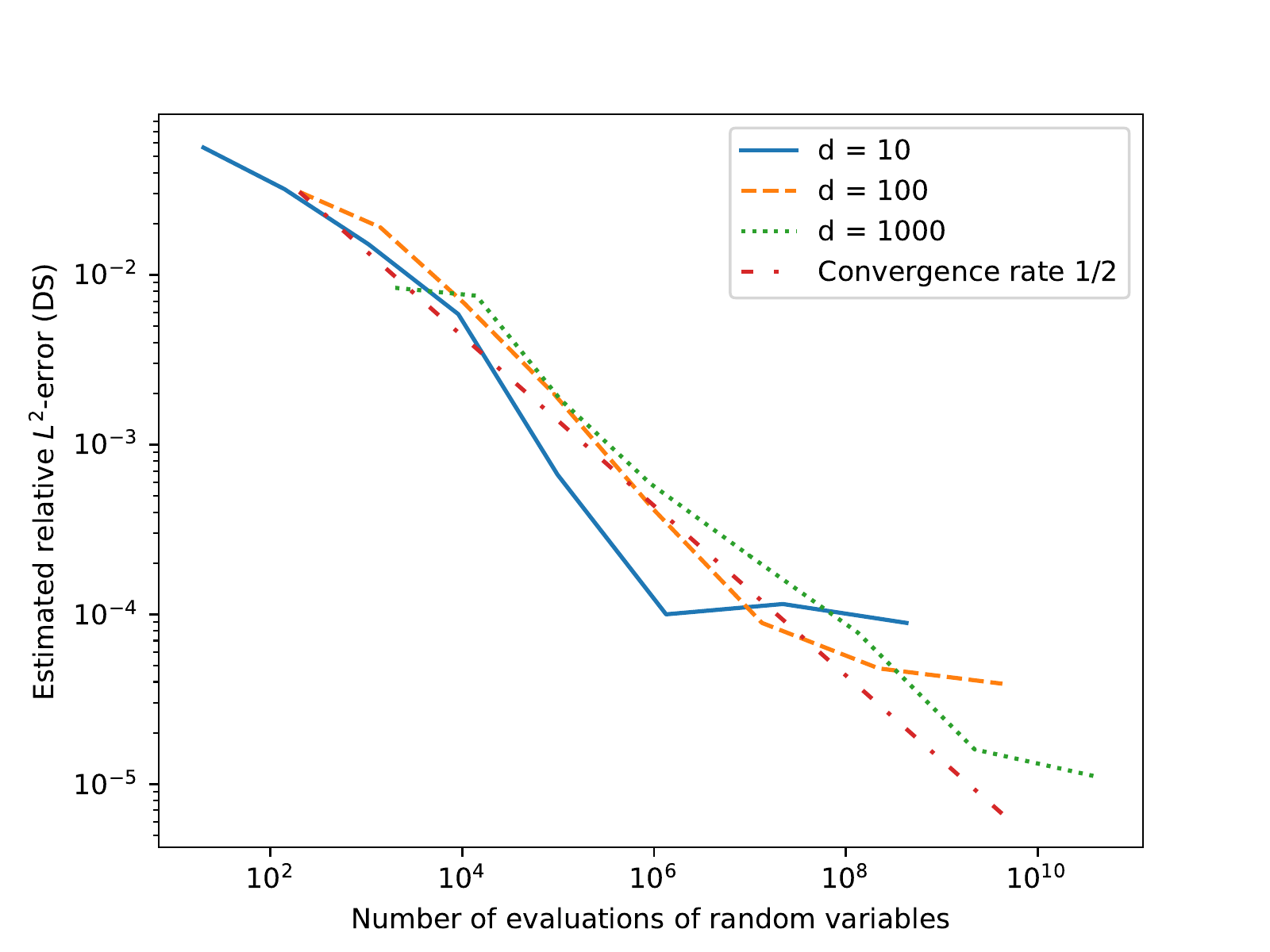} 
    \caption{Reference solutions computed by DS} 
    \label{fig:BlackScholes_DS} 
\end{subfigure} 
\begin{subfigure}{0.49\textwidth} 
    \includegraphics[width=\textwidth]{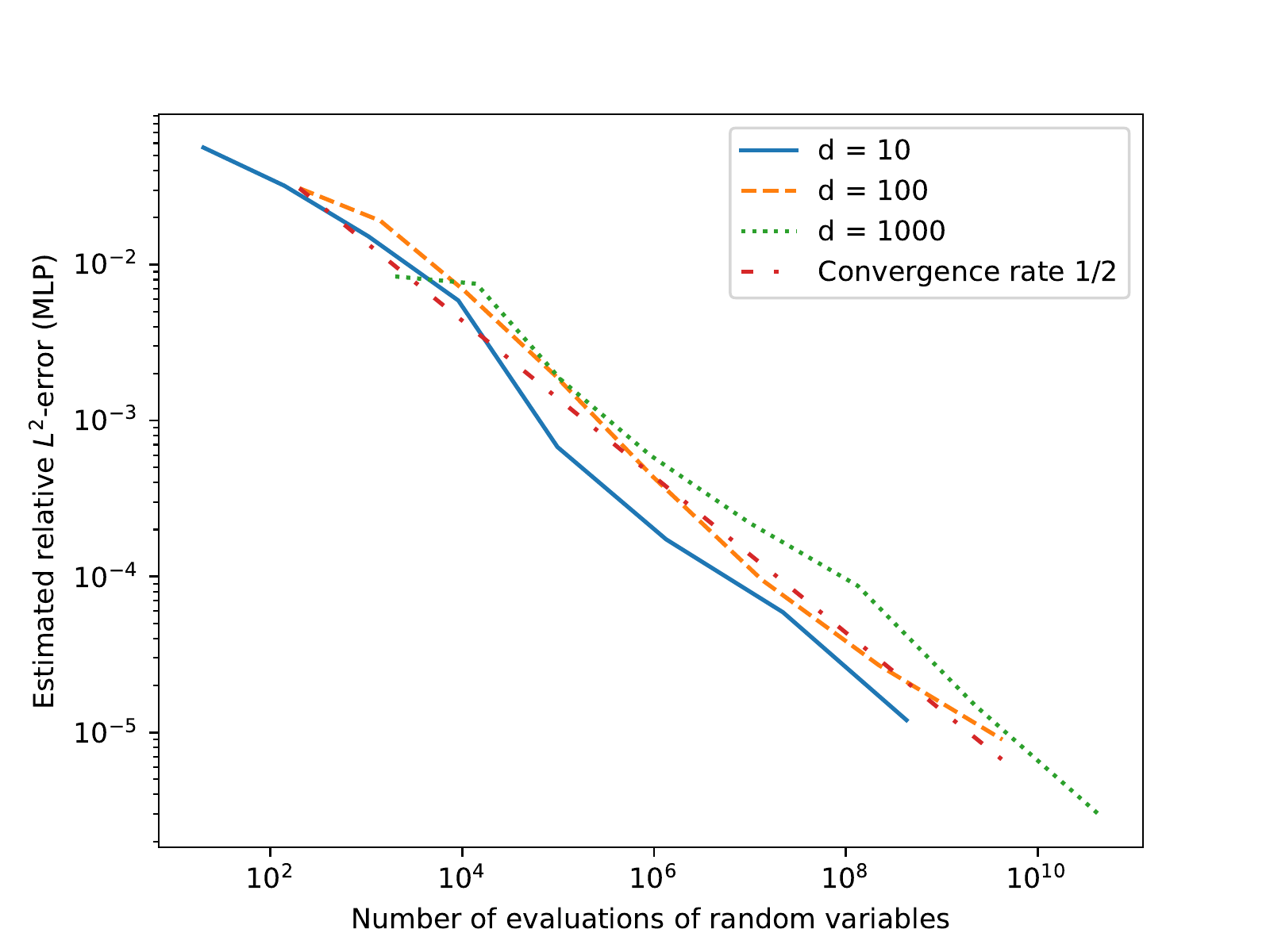} 
    \caption{Reference solutions computed by MLP}
    \label{fig:BlackScholes_MLP}
\end{subfigure}
\caption{Approximative plot of the relative $L^2$-error of the MLP approximation algorithm in \eqref{setting:MLP} against the computational effort of the algorithm in the case of the Black-Scholes PDE in \eqref{BlackScholes:eq1}.} 
\label{fig:BlackScholes} 
\end{figure}

Assume \cref{nogradsetting}, 
let $\xi = (50, \ldots, 50)\in \R^d$,
assume for all 
$t\in [0,T]$, 
$x=(x_{1},\dots,x_{d})\in \R^d$, 
$y\in \R$ 
that 
$k=1$, 
$T=1$, 
$f(x,y)= \tfrac{y}{1+y^{2}}$,  
$g(x) =\log(\tfrac 12 [1+\norm{x}^{2}])$, 
$\mu(x)=x$, 
and 
$\sigma(x) = \diag(x_{1},\dots,x_{d})$, 
let $\langle \cdot ,\cdot \rangle \colon (\cup_{q\in\N}(\R^{q}\times\R^{q}))\to \R$ be the standard scalar product,
and 
let $e_{1}, \dots, e_{d} \in \R^{d}$ 
satisfy that
$e_{1}= (1,0,\dots,0)$, \dots, $e_{d}=(0,\dots,0,1)$. 
Combining this, (\ref{setting:Xprocess}), and (\ref{setting:PDE}) with Hutzenthaler et al.\ \cite[Lemma 4.2]{VW2019defaultRisk} and the uniqueness property of solutions of stochastic differential equations (see, e.g., Klenke \cite{Klenke14})
ensures that for all 
$x=(x_{1},\dots,x_{d})\in \R^d$, 
$\theta \in \Theta$, 
$t\in [0,T]$, 
$s\in[t,T]$
it holds that 
\begin{equation}
\P\Big(X^{x,\theta}_{t,s} 
	= \big[x_{1} \exp \big( \tfrac {(s-t)}2 + \langle e_{1}, W^\theta_s - W^\theta_t\rangle \big), \dots, x_{d} \exp \big( \tfrac {(s-t)}2 + \langle e_{d}, W^\theta_s - W^\theta_t\rangle \big)\big] \Big)
	=1
\end{equation}
and 
\begin{equation}
\label{BlackScholes:eq1}
(\tfrac{\partial}{\partial t}u)(t,x) + f(x,u(t,x))+ \tfrac 12 \smallsum\limits^{d}_{i=1} \Big[ |x_{i}|^{2} \big(\tfrac{\del^{2}}{\del x_{i}\del x_{i}} u\big)(t,x)\Big]
+ \smallsum\limits_{i=1}^{d} \big[x_{i} (\tfrac{\partial}{\partial x_{i}} u)(t,x)\big] = 0.
\end{equation}
%
%
%%%%%%%%%%%%%%%%%%%%%%%%%%%%%%%%%%%%%%%%% 
%START 
%\input{01_Figure_descriptions/BlackScholes_description.tex}
%%%%%%%%%%%%%%%%%%%%%%%%%%%%%%%%%%%%%%%%% 
In \cref{BlackScholes-MLP} we approximately present for 
    $d \in \{10, 100, 1000 \}$, 
    $n \in \{1, 2, \ldots, 8 \}$ 
one random realization of $V^\theta_{n,n,\infty}(0,\xi)$ 
    (\nth {3} column in \cref{BlackScholes-MLP}), 
the relative $L^2$-error $\frac{( \E[ |V^{(0)}_{n,n,\infty}(0,\xi)-u(0,\xi)|^{2}])^{1/2}}{u(0,\xi)}$
    (\nth {5} and \nth {6} column in \cref{BlackScholes-MLP}),
the number of evaluations of one-dimensional random variables used to calculate one random realization of $V^\theta_{n,n,\infty}(0,\xi)$
    (\nth {7} column in \cref{BlackScholes-MLP}),
and
the runtime to calculate one random realization of $V^\theta_{n,n,\infty}(0,\xi)$
    (\nth {8} column in \cref{BlackScholes-MLP}).
In \cref{fig:BlackScholes} we approximately plot for
    $d \in \{10, 100, 1000 \}$,
     $n \in \{1, 2, \ldots, 8 \}$ 
the relative $L^2$-error $\frac{( \E[ |V^{(0)}_{n,n,\infty}(0,\xi)-u(0,\xi)|^{2}])^{1/2}}{u(0,\xi)}$ 
    (\nth {5} and \nth {6} column in \cref{BlackScholes-MLP})
against the number of evaluations of one-dimensional random variables used to calculate one random realization of $V^\theta_{n,n,\infty}(0,\xi)$
    (\nth {7} column in \cref{BlackScholes-MLP}).
The results in \cref{BlackScholes-MLP} and \cref{fig:BlackScholes} have been computed by means of \Cplusplus~code~\ref{code} in Section~\ref{sect:code} below.
For every 
    $n \in \{1, 2, \ldots, 8 \}$ 
for our approximative computations of the relative $L^2$-error $\frac{( \E[ |V^{(0)}_{n,n,\infty}(0,\xi)-u(0,\xi)|^{2}])^{1/2}}{u(0,\xi)}$ 
    (\nth {5} and \nth {6} column in \cref{BlackScholes-MLP})
the value $u(0,\xi)$ of the unknown exact solution 
in the relative $L^2$-error
has been approximated by means of an average of 5 independent runs of the deep splitting approximation method in Beck et al.\ \cite{beck2019deep}
    (\nth {5} column in \cref{BlackScholes-MLP})
and by means of an average of 5 independent evaluations of $V^{(0)}_{8,8,\infty}(0,\xi)$ 
    (\nth {6} column in \cref{BlackScholes-MLP}),
respectively,
and the expectation 
in the relative $L^2$-error 
has been approximated by means of Monte Carlo approximations involving 5 independent runs.

%%%%%%%%%%%%%%%%%%%%%%%%%%%%%%%%%%%%%%%%% 
%%END 
%\input{01_Figure_descriptions/BlackScholes_description.tex}
%%%%%%%%%%%%%%%%%%%%%%%%%%%%%%%%%%%%%%%%% 

%%%%%%%%%%%%%%%%%%%%%%%%%%%%%%%%%%%%%%%%%%%
%Numerical results
%%%%%%%%%%%%%%%%%%%%%%%%%%%%%%%%%%%%%%%%%%%

%%%%%%%%%%%%%%%%%%%%%%%%%%%%%%%%%%%%%%%%% 
%START 
%\input{../2_CSV_data/06_Simulations_200310/Merged_csv_files/Tabular/BlackScholes-MLP_tabular}
%%%%%%%%%%%%%%%%%%%%%%%%%%%%%%%%%%%%%%%%% 
\begin{table} 
\begin{center} 
\begin{tabular}{|c|c|c|c|c|c|c|c|}
\hline 
\thead{d} &\thead{n} &\thead{Result \\ of \\ MLP \\ algo- \\ rithm} &\thead{Refe- \\ rence \\ solu- \\ tions} &\thead{Esti- \\ mated \\ relative \\ $L^2$-error \\ (DS)} &\thead{Esti- \\ mated \\ relative \\ $L^2$-error \\ (MLP)} &\thead{Evaluations \\ of \\ random \\ variables} &\thead{Run- \\ time \\ in \\ sec- \\ onds} \\ 
\hline 
10 &1 &12.82440 &  &0.056260 &0.056265 &20 &0.00013 \\ 
10 &2 &11.95260 &DS: &0.032001 &0.031983 &140 &0.00012 \\ 
10 &3 &11.95160 &11.98841 &0.015161 &0.015158 &1050 &0.00051 \\ 
10 &4 &11.99310 &  &0.005870 &0.005885 &9080 &0.00329 \\ 
10 &5 &11.99990 &  &0.000666 &0.000676 &98300 &0.03670 \\ 
10 &6 &11.98910 &MLP: &0.000100 &0.000173 &1334340 &0.58832 \\ 
10 &7 &11.98600 &11.98736 &0.000115 &0.000059 &22032010 &4.92727 \\ 
10 &8 &11.98750 &  &0.000089 &0.000012 &428332080 &116.666 \\ 
\hline 
100 &1 &14.10430 &  &0.030842 &0.030856 &200 &0.00006 \\ 
100 &2 &14.75860 &DS: &0.019049 &0.019074 &1400 &0.00013 \\ 
100 &3 &14.52930 &14.68699 &0.006795 &0.006801 &10500 &0.00087 \\ 
100 &4 &14.66750 &  &0.001983 &0.001971 &90800 &0.00760 \\ 
100 &5 &14.68080 &  &0.000416 &0.000433 &983000 &0.09257 \\ 
100 &6 &14.68710 &MLP: &0.000089 &0.000095 &13343400 &1.53685 \\ 
100 &7 &14.68780 &14.68754 &0.000048 &0.000027 &220320100 &27.5609 \\ 
100 &8 &14.68760 &  &0.000039 &0.000009 &4283320800 &926.835 \\ 
\hline 
1000 &1 &16.92600 &  &0.008389 &0.008384 &2000 &0.00015 \\ 
1000 &2 &17.24340 &DS: &0.007540 &0.007536 &14000 &0.00067 \\ 
1000 &3 &17.11690 &17.07785 &0.001849 &0.001844 &105000 &0.00521 \\ 
1000 &4 &17.09480 &  &0.000591 &0.000600 &908000 &0.04907 \\ 
1000 &5 &17.07760 &  &0.000222 &0.000221 &9830000 &0.59165 \\ 
1000 &6 &17.07970 &MLP: &0.000078 &0.000087 &133434000 &8.48622 \\ 
1000 &7 &17.07750 &17.07766 &0.000016 &0.000015 &2203201000 &231.843 \\ 
1000 &8 &17.07770 &  &0.000011 &0.000003 &42833208000 &8795.39 \\ 
\hline 
\end{tabular} 
\end{center} 
\caption{ \label{BlackScholes-MLP} Numerical simulations for the MLP approximation algorithm in \eqref{setting:MLP} in the case of the Black-Scholes PDE in \eqref{BlackScholes:eq1}} 
\end{table} 

%%%%%%%%%%%%%%%%%%%%%%%%%%%%%%%%%%%%%%%%% 
%%END 
%\input{../2_CSV_data/06_Simulations_200310/Merged_csv_files/Tabular/BlackScholes-MLP_tabular}
%%%%%%%%%%%%%%%%%%%%%%%%%%%%%%%%%%%%%%%%% 

%%%%%%%%%%%%%%%%%%%%%%%%%%%%%%%%%%%%%%%%% 
%START 
%\input{../2_CSV_data/06_Simulations_200310/Merged_csv_files/Plots/BlackScholes_plotfigure}
%%%%%%%%%%%%%%%%%%%%%%%%%%%%%%%%%%%%%%%%% 

%%%%%%%%%%%%%%%%%%%%%%%%%%%%%%%%%%%%%%%%% 
%%END 
%\input{../2_CSV_data/06_Simulations_200310/Merged_csv_files/Plots/BlackScholes_plotfigure}
%%%%%%%%%%%%%%%%%%%%%%%%%%%%%%%%%%%%%%%%% 

\begin{comment}
%%%%%%%%%%%%%%%%%%%%%%%%%%%%%%%%%%%%%%%%%%%
\subsection{Hamilton-Jacobi-Bellman equation (under construction)}
%%%%%%%%%%%%%%%%%%%%%%%%%%%%%%%%%%%%%%%%%%%

In this subsection we consider the Hamilton-Jacobi-Bellman equation (cf.\ Bellman \cite{Bellman1957})

Assume Setting [+++]  
and 
assume for all 
$t\in [0,T]$, 
$x,v\in \R^d$, 
$y\in \R$, 
$z\in\R^{d}$ 
that 
$d = 100$, 
$k=1$, 
$T=1$, 
$f(t,x,y,z)= -\norm{z}^{2}$, 
$g(x) =\log(\tfrac 12 [1+\norm{x}^{2}])$, 
$\mu(x)=0$, 
and 
$\sigma(x)v = v$. 
%
This, (\ref{setting:Xprocess}), and (\ref{setting:PDE}) 
ensure that 
for all 
$x\in \R^d$, 
$\theta \in \Theta$, 
$t\in [0,T]$, 
$s\in[t,T]$ 
it holds that 
$\P\big(X^{x,\theta}_{t,s} = x+W^{\theta}_{s}-W^{\theta}_{t}\big)=1$ 
and
\begin{equation}
(\tfrac{\partial}{\partial t}u)(t,x) + (\Delta_{x}u)(t,x) + f\big(t,x,u(t,x),(\nabla_x u)(t,x)\big) = 0.
\end{equation}
%
Consider the function
\begin{equation}
 \{1,\dots,7\}\ni n \mapsto \Big( \E\Big[ \big|V^{(0)}_{n,n,n}(0,0)-u(0,0)\big|^{2}\Big]\Big)^{\!\nicefrac 12} \in \R.
\end{equation}

%%%%%%%%%%%%%%%%%%%%%%%%%%%%%%%%%%%%%%%%%%%
%Numerical results
%%%%%%%%%%%%%%%%%%%%%%%%%%%%%%%%%%%%%%%%%%%
%\csvautotabular{../2_CSV_data/1_First_Simulations/Sine-Gordon_equation.csv}

\end{comment}

%%%%%%%%%%%%%%%%%%%%%%%%%%%%%%%%%%%%%%%%%%%
% Source Codes
%%%%%%%%%%%%%%%%%%%%%%%%%%%%%%%%%%%%%%%%%%%

\section{Source code}
\label{sect:code}

In this section we present the source code (see \Cplusplus~code~\ref{code} below) which was employed to produce the results in \cref{subsect:AllenCahn,subsect:SineGordon,subsect:PDESystem,subsect:BlackScholes} above.
All of the numerical simulations presented in \cref{subsect:AllenCahn,subsect:SineGordon,subsect:PDESystem,subsect:BlackScholes} above were built and run on a system with an AMD Ryzen 9 3950X 16c/32t and 64 GB DDR4-3600 memory running Ubuntu 19.10. The provided source code uses the Eigen \Cplusplus~Library (version 3.3.7) and the POSIX Threads API to allow for parallelism on modern multicore CPUs. It was compiled with the \Cplusplus~compiler of the GNU Compiler Collection (version 7.5.0) with optimization level 3 (\verb+-O3+). The different examples can be selected at compile time by providing a preprocessor symbol using the \verb+-D+ option to activate the corresponding preprocessor macro. 
Possible choices for the preprocessor symbol are 
\verb+ALLEN_CAHN+ (see \cref{subsect:AllenCahn}), 
\verb+SINE_GORDON+ (see \cref{subsect:SineGordon}), 
\verb+PDE_SYSTEM+ (see \cref{subsect:PDESystem}), and 
\verb+SEMILINEAR_BS+ (see \cref{subsect:BlackScholes}). 
For example, the source code for the Allen-Cahn example (see \cref{subsect:AllenCahn}) was compiled using the command:
%\vspace{-0.1cm}
%\begin{verbatim}
%g++ -DALLEN_CAHN -O3 -o mlp mlp.cpp -lpthread
%\end{verbatim}
%\vspace{-0.1cm}
\verb{g++ -DALLEN_CAHN -O3 -o mlp mlp.cpp -lpthread{.
%TODO: ask sebastian if the above formulation is ok
Note that if the Eigen headers are not available system-wide the path has to be provided using the \verb+-I+ option in the command above.

\begin{lstlisting}[ label = code, caption = Source code for \cref{subsect:AllenCahn,subsect:SineGordon,subsect:PDESystem,subsect:BlackScholes}]
#define N_MAX 8

#ifdef ALLEN_CAHN
#define eq_name "Allen-Cahn_equation"
#define rdim 1
#define g(x) ArrayXd tmp = ArrayXd::Zero(1, 1); tmp(0) = 1. / (2. + 2. / 5. * x.square().sum())
#define X_sde(s, t, x, w) x + sqrt(2. * (t - s)) * w
#define fn(y) ArrayXd ret = ArrayXd::Zero(1, 1); double phi_r = std::min(4., std::max(-4., y(0))); ret(0) = phi_r - phi_r * phi_r * phi_r
#endif
#ifdef SINE_GORDON
#define eq_name "Sine-Gordon_equation"
#define rdim 1
#define g(x) ArrayXd tmp = ArrayXd::Zero(1, 1); tmp(0) = 1. / (2. + 2. / 5. * x.square().sum())
#define X_sde(s, t, x, w) x + sqrt(2. * (t - s)) * w
#define fn(y) ArrayXd ret = ArrayXd::Zero(1, 1); ret(0) = sin(y(0))
#endif
#ifdef SEMILINEAR_BS
#define eq_name "Semilinear_Black-Scholes_equation"
#define rdim 1
#define g(x) ArrayXd tmp = ArrayXd::Zero(1, 1); tmp(0) = log(0.5 * (1. + x.square().sum()))
#define X_sde(s, t, x, w) x * ((t - s) / 2. + sqrt(t - s) * w).exp()
#define fn(y) ArrayXd ret = ArrayXd::Zero(1, 1); ret(0) = y(0) / (1. + y(0) * y(0))
#endif
#ifdef PDE_SYSTEM
#define eq_name "Semilinear_PDE_system"
#define rdim 2
#define X_sde(s, t, x, w) x + sqrt(2. * (t - s)) * w
#define g(x) ArrayXd tmp = ArrayXd::Zero(2, 1); tmp(0) = 1. / (2. + 2. / 5. * x.square().sum()); tmp(1) = log(0.5 * (1. + x.square().sum()))
#define fn(y) ArrayXd ret = ArrayXd::Zero(2, 1); ret(0) = y(1) / (1. + y(1) * y(1)); ret(1) = 2. * y(0) / 3.
#endif


#include <iostream>
#include <iomanip>
#include <fstream>
#include <cstdint>
#include <random>
#include <cmath>
#include <ctime>
#include <thread>
#include <chrono>
#include <eigen3/Eigen/Dense>

using Eigen::ArrayXd;

struct mlp_t {
	uint8_t m;
	uint8_t n;
	uint8_t l;
	uint16_t d;
	ArrayXd x;
	double s;
	double t;
	ArrayXd res;
};


ArrayXd f(const ArrayXd &v);
ArrayXd mlp_call(uint8_t m, uint8_t n, uint8_t l, uint16_t d, ArrayXd &x, double s, double t);
ArrayXd ml_picard(uint8_t m, uint8_t n, uint16_t d, ArrayXd &x, double s, double t, bool start_threads);
void mlp_thread(mlp_t &mlp_args);

int main() {

	std::string s = eq_name;
	std::cout << s << std::endl << std::endl << std::setprecision(8);

	std::ofstream out_file;
	out_file.open(s + "_mlp.csv");
	out_file << "d, T, n, ";
	for (uint8_t i=0; i < rdim; i++) {
		out_file << "result_" << (int)i << ", ";
	}
	out_file << "elapsed_secs" << std::endl;

	double T[1] = {1.};
	uint16_t d[3] = {10, 100, 1000};

	for (uint16_t j = 0; j < sizeof(d) / sizeof(d[0]); j++) {
		for (uint8_t k = 0; k < sizeof(T) / sizeof(T[0]); k++) {
			for (uint8_t n = 1; n <= N_MAX; n++) {
				std::chrono::time_point<std::chrono::high_resolution_clock> start_time = std::chrono::high_resolution_clock::now();
		#if defined(NONLINEAR_BS) || defined(SEMILINEAR_BS)
				ArrayXd xi = ArrayXd::Constant(d[j], 1, 50.);
		#else
				ArrayXd xi = ArrayXd::Zero(d[j], 1);
		#endif
				ArrayXd result = ml_picard(n, n, d[j], xi, 0., T[k], true);

				std::chrono::time_point<std::chrono::high_resolution_clock> end_time = std::chrono::high_resolution_clock::now();
				double elapsed_secs = double(std::chrono::duration_cast<std::chrono::microseconds>(end_time - start_time).count()) / 1000. / 1000.;
				std::cout << "T: " << T[k] << std::endl << "d: " << (int)d[j] << std::endl;
				std::cout << "n: " << (int)n << std::endl << "Result:" << std::endl << result << std::endl;
				std::cout << "Elapsed secs: " << elapsed_secs << std::endl << std::endl;

				out_file << (int)d[j] << ", " << T[k] << ", " << (int)n << ", ";
				for (uint8_t i = 0; i < rdim; i++) {
					out_file << result(i) << ", ";
				}
				out_file << elapsed_secs << std::endl;
			}
		}
	}

	out_file.close();

	return 0;
}

ArrayXd f(const ArrayXd &v) {
	fn(v);
	return ret;
}

void mlp_thread(mlp_t &mlp_args) {
	mlp_args.res = mlp_call(mlp_args.m, mlp_args.n, mlp_args.l, mlp_args.d, mlp_args.x, mlp_args.s, mlp_args.t);
}

ArrayXd mlp_call(uint8_t m, uint8_t n, uint8_t l, uint16_t d, ArrayXd &x, double s, double t) {
	ArrayXd a = ArrayXd::Zero(rdim, 1);
	ArrayXd b = ArrayXd::Zero(rdim, 1);
	double r = 0.;
	ArrayXd x2;
	uint32_t num;
	static thread_local std::mt19937 generator(128 + clock() + std::hash<std::thread::id>()(std::this_thread::get_id()));
	static thread_local std::normal_distribution<> normal_distribution{0., 1.};
	static thread_local std::uniform_real_distribution<double> uniform_distribution(0., 1.);
	if (l < 2) {
		num = (uint32_t)(pow(m, n - l) + 0.5);
		for (uint32_t k = 0; k < num; k++) {
			r = s + (t - s) * uniform_distribution(generator);
			x2 = ArrayXd::NullaryExpr(d, [&](){ return normal_distribution(generator); });
			x2 = X_sde(s, r, x, x2);
			b += f(ml_picard(m, l, d, x2, r, t, false));
		}
		a += (t - s) * (b / ((double)num));
	} else {
		num = (uint32_t)(pow(m, n - l) + 0.5);
		for (uint32_t k = 0; k < num; k++) {
			r = s + (t - s) * uniform_distribution(generator);
			x2 = ArrayXd::NullaryExpr(d, [&](){ return normal_distribution(generator); });
			x2 = X_sde(s, r, x, x2);
			b += (f(ml_picard(m, l, d, x2, r, t, l > m - 5)) - f(ml_picard(m, l - 1, d, x2, r, t, l > m - 5)));
		}
		a += (t - s) * (b / ((double)num));
	}
	return a;
}

ArrayXd ml_picard(uint8_t m, uint8_t n, uint16_t d, ArrayXd &x, double s, double t, bool start_threads) {

	if (n == 0) return ArrayXd::Zero(rdim, 1);

	ArrayXd a = ArrayXd::Zero(rdim, 1);
	ArrayXd a2 = ArrayXd::Zero(rdim);
	ArrayXd b = ArrayXd::Zero(rdim, 1);

	double r = 0.;
	std::thread threads[16];
	mlp_t mlp_args[16];
	ArrayXd x2;
	uint32_t num;
	static thread_local std::mt19937 generator(clock() + std::hash<std::thread::id>()(std::this_thread::get_id()));
	static thread_local std::normal_distribution<> normal_distribution{0., 1.};
	static thread_local std::uniform_real_distribution<double> uniform_distribution(0., 1.);


	if (start_threads) {

		for (uint8_t l = 0; l < n; l++) {
			mlp_t mlp_arg;
			mlp_arg.m = m;
			mlp_arg.n = n;
			mlp_arg.l = l;
			mlp_arg.d = d;
			mlp_arg.x = x.replicate(1, 1);
			mlp_arg.s = s;
			mlp_arg.t = t;
			mlp_arg.res = 0.;
			mlp_args[l] = mlp_arg;
			threads[l] = std::thread(mlp_thread, std::ref(mlp_args[l]));
		}

		num = (uint32_t)(pow(m, n) + 0.5);
		for (uint32_t k = 0; k < num; k++) {
			x2 = ArrayXd::NullaryExpr(d, [&](){ return normal_distribution(generator); });
			x2 = X_sde(s, t, x, x2);
			g(x2);
			a2 += tmp;
		}

		a2 /= (double)num;

		for (uint8_t l = 0; l < n; l++) {
			threads[l].join();
			a += mlp_args[l].res;
		}

	} else {

		for (uint8_t l = 0; l < std::min(n, (uint8_t)2); l++) {
			b = ArrayXd::Zero(rdim, 1);
			num = (uint32_t)(pow(m, n - l) + 0.5);
			for (uint32_t k = 0; k < num; k++) {
				r = s + (t - s) * uniform_distribution(generator);
				x2 = ArrayXd::NullaryExpr(d, [&](){ return normal_distribution(generator); });
				x2 = X_sde(s, r, x, x2);
				b += f(ml_picard(m, l, d, x2, r, t, false));
			}
			a += (t - s) * (b / ((double)num));
		}

		for (uint8_t l = 2; l < n; l++) {
			b = ArrayXd::Zero(rdim, 1);
			num = (uint32_t)(pow(m, n - l) + 0.5);
			for (uint32_t k = 0; k < num; k++) {
				r = s + (t - s) * uniform_distribution(generator);
				x2 = ArrayXd::NullaryExpr(d, [&](){ return normal_distribution(generator); });
				x2 = X_sde(s, r, x, x2);
				b += (f(ml_picard(m, l, d, x2, r, t, false)) - f(ml_picard(m, l - 1, d, x2, r, t, false)));
			}
			a += (t - s) * (b / ((double)num));
		}

		num = (uint32_t)(pow(m, n) + 0.5);
		for (uint32_t k = 0; k < num; k++) {
			x2 = ArrayXd::NullaryExpr(d, [&](){ return normal_distribution(generator); });
			x2 = X_sde(s, t, x, x2);
			g(x2);
			a2 += tmp;
		}

		a2 /= (double)num;

	}

	return a + a2;
}
\end{lstlisting}

%%%%%%%%%%%%%%%%%%%%%%%%%%%%%%%%%%%%%%%%%%%
%%%%%%%%%%%%%%%%%%%%%%%%%%%%%%%%%%%%%%%%%%%
%Bibliography
%%%%%%%%%%%%%%%%%%%%%%%%%%%%%%%%%%%%%%%%%%%
%%%%%%%%%%%%%%%%%%%%%%%%%%%%%%%%%%%%%%%%%%%

\bibliographystyle{acm}
%%%%%%%%%%%%%%%%%%%%%%%%%%%%%%%%%%%%%%%%% 
%\bibliography{../../0A_Templates_and_supporting_material/3_Bibfile/Main_bibfile.bib}

%%%%%%%%%%%%%%%%%%%%%%%%%%%%%%%%%%%%%%%%% 

\end{document}